\documentclass{amsart}
\usepackage{amssymb}
\usepackage{amscd}
\usepackage{verbatim}
\usepackage{epsfig}
\begin{document}
\newcommand\Mand{\ \text{and}\ }
\newcommand\Mwith{\ \text{with}\ }
\newcommand\Mfor{\ \text{for}\ }
\newcommand\Mst{\ \text{such that}\ }
\newcommand\Mor{\ \text{or}\ }
\newcommand\Mif{\ \text{if}\ }
\newcommand\Miff{\ \text{iff}\ }
\newcommand\Mthen{\ \text{then}\ }
\newcommand\nin{\notin}
\newcommand\identity{\operatorname{id}}
\newcommand\Id{\operatorname{Id}}
\newcommand\Real{\mathbb{R}}
\newcommand\RR{\mathbb{R}}
\newcommand\pos{\Real^+}
\newcommand\Rnp{\Real\setminus\{0\}}
\newcommand\nzero{\setminus\{0\}}
\newcommand\Cx{\mathbb{C}}
\newcommand\Cxp{\Cx^+}
\newcommand\Cxm{\Cx^-}
\newcommand\Nat{\mathbb{N}}
\newcommand\halfNat{{\frac{1}{2}}\mathbb{N}}
\newcommand\intgr{\mathbb{Z}}
\newcommand\HH{\mathbb{H}}
\newcommand\im{\operatorname{Im}}
\newcommand\re{\operatorname{Re}}
\newcommand\sign{\operatorname{sign}}
\newcommand\codim{\operatorname{codim}}
\newcommand\End{\operatorname{End}}
\newcommand\Ker{\operatorname{Ker}}
\newcommand\Hom{\operatorname{Hom}}
\newcommand\tr{\operatorname{tr}}
\newcommand\Tr{\operatorname{Tr}}
\newcommand\ideal{{\mathcal I}}
\newcommand\Span{\operatorname{span}}
\newcommand\image{\operatorname{image}}
\newcommand\Range{\operatorname{Ran}}
\newcommand\Graph{\operatorname{graph}}
\newcommand\slim{\operatornamewithlimits{s-lim}}
\newcommand\spp{\operatorname{sp}}
\newcommand\sll{\operatorname{sl}}
\newcommand\sol{\operatorname{so}}
\newcommand\SL{\operatorname{SL}}
\newcommand\SO{\operatorname{SO}}
\newcommand\On{\operatorname{O}}
\newcommand\pa{\partial}
\newcommand\eff{\mathrm{eff}}
\newcommand\Rn{\Real^n}
\newcommand\Rm{\Real^m}
\newcommand\RN{\Real^N}
\newcommand\RtN{\Real^{2N}}
\newcommand\RM{\Real^M}
\newcommand\sphere{\mathbb{S}}
\newcommand\Sn{\sphere^{n-1}}
\newcommand\Sm{\sphere^{m-1}}
\newcommand\Snp{\sphere^n_+}
\newcommand\Smp{\sphere^m_+}
\newcommand\SN{\sphere^{N-1}}
\newcommand\SNp{\sphere^N_+}
\newcommand\circlep{\sphere^1_+}
\newcommand\Phom{P_{h}}
\newcommand\Shom{S_{h}}
\newcommand\distance{\operatorname{dist}}
\newcommand\cl{\operatorname{cl}}
\newcommand\interior{\operatorname{int}}
\newcommand\Fa{\operatorname{Fa}}
\newcommand\ff{\operatorname{ff}}
\newcommand\mf{\operatorname{mf}}
\newcommand\cf{\operatorname{cf}}
\newcommand\scf{\operatorname{sf}}
\newcommand\lf{\operatorname{lf}}
\newcommand\rf{\operatorname{rf}}
\newcommand\indfam{{\mathcal K}}
\newcommand\fraka{{\mathfrak a}}
\newcommand\calA{{\mathcal A}}
\newcommand\calB{{\mathcal B}}
\newcommand\calR{{\mathcal R}}
\newcommand\calO{{\mathcal O}}
\newcommand\calJ{{\mathcal J}}
\newcommand\calM{{\mathcal M}}
\newcommand\calN{{\mathcal N}}
\newcommand\calX{{\mathcal X}}
\newcommand\calU{{\mathcal U}}
\newcommand\calV{{\mathcal V}}
\newcommand\calF{{\mathcal F}}
\newcommand\calG{{\mathcal G}}
\newcommand\calT{{\mathcal T}}
\newcommand\calC{{\mathcal C}}
\newcommand\calCt{{\tilde {\mathcal C}}}
\newcommand\calCL{{\mathcal C}_{\text L}}
\newcommand\calCR{{\mathcal C}_{\text R}}
\newcommand\Cinf{{\mathcal C}^{\infty}}
\newcommand\dist{{\mathcal C}^{-\infty}}
\newcommand\dCinf{{\dot{\mathcal C}}^{\infty}}
\newcommand\ddist{\dot\dist}
\newcommand\Cj{{\mathcal C}^j}
\newcommand\Linf{L^{\infty}}
\newcommand\phg{{\text{phg}}}
\newcommand\comp{{\text{comp}}}
\newcommand\bcon{{\mathcal A}}
\newcommand\bconc{{\mathcal A}_{\text{phg}}}
\newcommand\Sch{{\mathcal S}}
\newcommand\temp{\Sch^{\prime}}
\newcommand\Diff{\operatorname{Diff}}
\newcommand\Diffb{\operatorname{Diff}_{\text{b}}}
\newcommand\Diffc{\operatorname{Diff}_{\text{c}}}
\newcommand\Diffsc{\operatorname{Diff}_{\text{sc}}}
\newcommand\DiffI{\operatorname{Diff}_{\text{I}}}
\newcommand\DiffIq{\operatorname{Diff}_{\text{I},q}}
\newcommand\sing{\text{sing}}
\newcommand\reg{\text{reg}}
\newcommand\supp{\operatorname{supp}}
\newcommand\ssupp{\operatorname{sing\ supp}}
\newcommand\csupp{\operatorname{cone\ supp}}
\newcommand\esupp{\operatorname{ess\ supp}}
\newcommand\Fr{{\mathcal F}}
\newcommand\Frinv{\Fr^{-1}}
\newcommand\bop{{\mathcal B}}
\newcommand\spec{\operatorname{spec}}
\newcommand\pspec{\spec_{pp}}
\newcommand\cspec{\spec_{c}}
\newcommand\FIO{{\mathcal I}}
\newcommand\SP{\operatorname{RC}}
\newcommand\RC{\operatorname{RC}}
\newcommand\Symc{S_c}
\newcommand\Symca{S_c^{\alpha}}
\newcommand\Symczero{S_c^{0,...,0}}
\newcommand\sci{{}^{\text{sc}}}
\newcommand\sct{\sci T^*}
\newcommand\scT{\sci T}
\newcommand\scdt{\sci \dot T^*}
\newcommand\dS{\dot S^*}
\newcommand\dT{\dot T^*}
\newcommand\dSreg{\dot\Sigma_{\text reg}}
\newcommand\scct{\sci\bar{T}^*}
\newcommand\Csc{C_{\text{sc}}}
\newcommand\SNpscd{(\SNp)^2_{\text{sc}}}
\newcommand\scdiag{\Delta_{\text{sc}}}
\newcommand\projscl{\pi^L_{\text{sc}}}
\newcommand\projscr{\pi^R_{\text{sc}}}
\newcommand\scHL{\sci H^{2,0}_{|\zeta|^2-\lambda^2}}
\newcommand\scHrg{\sci H^{2,0}_{\sqrt{g}}}
\newcommand\Hsc{H_{\text{sc}}}
\newcommand\Char{\operatorname{Char}}
\newcommand\dChar{\operatorname{\dot Char}}
\newcommand\WF{\operatorname{WF}}
\newcommand\WFp{\operatorname{WF^{\prime}}}
\newcommand\WFsc{\operatorname{WF}_{\text{sc}}}
\newcommand\WFscp{\operatorname{WF_{sc}^{\prime}}}
\newcommand\WFC{\operatorname{WF}_C}
\newcommand\WFCi{\operatorname{WF}_{C_i}}
\newcommand\elliptic{\operatorname{ell}}
\newcommand\Psop{\operatorname{\Psi}}
\newcommand\Psiscrs{\operatorname{\Psi_{sc}^{-2,\infty}}}
\newcommand\Psiscr{\operatorname{\Psi_{sc}^{-2,0}}}
\newcommand\Psiscrm{\operatorname{\Psi_{sc}^{0,2}}}
\newcommand\PsiscHam{\operatorname{\Psi_{sc}^{2,0}}}
\newcommand\Psisci{\operatorname{\Psi_{sc}^{*,*}}}
\newcommand\Psiscid{\operatorname{\Psi_{sc}^{0,0}}}
\newcommand\Psiscis{\operatorname{\Psi_{sc}^{0,\infty}}}
\newcommand\Psiscsi{\operatorname{\Psi_{sc}^{-\infty,0}}}
\newcommand\Psiscs{\operatorname{\Psi_{sc}^{-\infty,\infty}}}
\newcommand\Psiscalg{\operatorname{\Psi_{sc}^{\infty,-\infty}}}
\newcommand\nullHam{{\mathcal N}}
\newcommand\charD{\Sigma_{\Delta-\lambda^2}}
\newcommand\charLap{\Sigma_{\Delta-\lambda}}
\newcommand\Snl{\Sn_{\lambda}}
\newcommand\SNl{\SN_{\lambda}}
\newcommand\gammat{\tilde\gamma}
\newcommand\gammasc{\gamma}
\newcommand\Tau{\mathcal{T}}
\newcommand\taut{\tilde\tau}
\newcommand\taub{\bar\tau}
\newcommand\Nout{N^+_{\lambda}}
\newcommand\Nin{N^-_{\lambda}}
\newcommand\Nio{N^{\pm}_{\lambda}}
\newcommand\El{E_{\lambda}}
\newcommand\Elt{\tilde E_{\lambda}}
\newcommand\Eil{E^i_{\lambda}}
\newcommand\Ejl{E^j_{\lambda}}
\newcommand\Eajl{E^{\alpha_j}_{\lambda}}
\newcommand\Eilt{\tilde E^i_{\lambda}}
\newcommand\Np{N^+}
\newcommand\Nm{N^-}
\newcommand\Npm{N^{\pm}}
\newcommand\Fin{F^-(\lambda)}
\newcommand\Fini{F^-_i(\lambda)}
\newcommand\Fout{F^+(\lambda)}
\newcommand\Fouti{F^+_i(\lambda)}
\newcommand\Foutj{F^+_j(\lambda)}
\newcommand\Rout{R^+_{\lambda}}
\newcommand\Routl{R^+_{\lambda^2}}
\newcommand\Routsgnl{R^{\sign\lambda}_{\lambda^2}}
\newcommand\Rin{R^-_{\lambda}}
\newcommand\Rinl{R^-_{\lambda^2}}
\newcommand\Rinsgnl{R^{-\sign\lambda}_{\lambda^2}}
\newcommand\Rio{R^{\pm}_{\lambda}}
\newcommand\Riol{R^{\pm}_{\lambda^2}}
\newcommand\Roi{R^{\mp}_{\lambda}}
\newcommand\Roil{R^{\mp}_{\lambda^2}}
\newcommand\Riob{R^{\pm}}
\newcommand\Roib{R^{\mp}}
\newcommand\Tio{T^{\pm}}
\newcommand\Tiob{T^{\pm}_{\ff}}
\newcommand\Toi{T^{\mp}}
\newcommand\Toib{T^{\mp}_{\ff}}
\newcommand\TIiob{T_I^{\pm}}
\newcommand\Rinb{R^-}
\newcommand\Rinbsgnl{R^{-\sign\lambda}}
\newcommand\Tin{T^-}
\newcommand\Tinb{T^-_{\ff}}
\newcommand\TIinb{T^-_I}
\newcommand\Routb{R^+}
\newcommand\Routbsgnl{R^{\sign\lambda}}
\newcommand\Tout{T^+}
\newcommand\Toutb{T^+_{\ff}}
\newcommand\TIoutb{T^+_I}
\newcommand\Rlkf{(|\xib|^2-(\lambda-i0)^2)^{-1}}
\newcommand\Rlk{\rho_0(\lambda)}
\newcommand\Rmlk{\rho_0(-\lambda)}
\newcommand\Rpmlk{\rho_0(\pm\lambda)}
\newcommand\Rlka{\rho_1(\lambda)}
\newcommand\Rlkb{\rho_2(\lambda)}
\newcommand\Rilk{\rho_i(\lambda)}
\newcommand\reduced{\natural}
\newcommand\Rlf{R_0(\lambda)}
\newcommand\Rla{R_1(\lambda)}
\newcommand\Rlb{R_2(\lambda)}
\newcommand\Ril{R_i(\lambda)}
\newcommand\Rlj{R_j(\lambda)}
\newcommand\Rlft{R_0(\lambda)}
\newcommand\Rflambda{R_0^{\reduced}(\sigma)}
\newcommand\RV{R^{\reduced}_V}
\newcommand\Rfsigma{R_0^{\reduced}(\sigma)}
\newcommand\Rfsigmah{R_0^{\reduced}(\sigma^{1/2})}
\newcommand\Rfzero{R_0^{\reduced}(0)}
\newcommand\RlV{R^{\reduced}_V(\sigma)}
\newcommand\RlVi{R^{\reduced}_{V_i}(\sigma)}
\newcommand\RlVt{R_V(\lambda)}
\newcommand\RlVtL{{R}_V^L(\lambda)}
\newcommand\RlVtR{{R}_V^R(\lambda)}
\newcommand\RlVit{{R}_{V_i}(\lambda)}
\newcommand\RlVta{{R}_V^{(1)}(\lambda)}
\newcommand\RlVtk{{R}_V^{(k)}(\lambda)}
\newcommand\RlVatV{{R}_{V_{\alpha}}(\lambda)V_{\alpha}}
\newcommand\RlVatVa{{R}_{V_{\alpha_1}}(\lambda)V_{\alpha_1}}
\newcommand\RlVatVb{{R}_{V_{\alpha_2}}(\lambda)V_{\alpha_2}}
\newcommand\RlVatVk{{R}_{V_{\alpha_k}}(\lambda)V_{\alpha_k}}
\newcommand\RlVatVkk{{R}_{V_{\alpha_{k+1}}}(\lambda)V_{\alpha_{k+1}}}
\newcommand\RlVaptV{{R}_{V_{\alpha'}}(\lambda)V_{\alpha'}}
\newcommand\RlVapptV{{R}_{V_{\alpha''}}(\lambda)V_{\alpha''}}
\newcommand\RlVajtV{{R}_{V_{\alpha_j}}(\lambda)V_{\alpha_j}}
\newcommand\RlVaktV{{R}_{V_{\alpha_k}}(\lambda)V_{\alpha_k}}
\newcommand\RlVakktV{{R}_{V_{\alpha_{k+1}}}(\lambda)V_{\alpha_{k+1}}}
\newcommand\Tl{T(\lambda)}
\newcommand\Tlt{\tilde\Tl}
\newcommand\Tltp{\tilde T'(\lambda)}
\newcommand\Tltpp{\tilde T''(\lambda)}
\newcommand\Tli{T_i(\lambda)}
\newcommand\Tlit{\tilde\Tli}
\newcommand\Tlip{T_i'(\lambda)}
\newcommand\Tlipp{T_i''(\lambda)}
\newcommand\Tlj{T_j(\lambda)}
\newcommand\Tla{T_{\alpha}(\lambda)}
\newcommand\Tlaa{T_{\alpha_1}(\lambda)}
\newcommand\Tlab{T_{\alpha_2}(\lambda)}
\newcommand\Tlak{T_{\alpha_k}(\lambda)}
\newcommand\Tlakt{\tilde\Tlak}
\newcommand\Tlaj{T_{\alpha_j}(\lambda)}
\newcommand\Tlajj{T_{\alpha_{j+1}}(\lambda)}
\newcommand\Tlajp{T_{\alpha_j}'(\lambda)}
\newcommand\Tlajpt{\tilde\Tlajp}
\newcommand\Tlajt{\tilde\Tlaj}
\newcommand\Tlakk{T_{\alpha_{k+1}}(\lambda)}
\newcommand\Tlakkp{T_{\alpha_{k+1}}'(\lambda)}
\newcommand\Tlap{T_{\alpha'}(\lambda)}
\newcommand\Tlapt{\tilde\Tlap}
\newcommand\Tlapp{T_{\alpha''}(\lambda)}
\newcommand\Tkl{T^{(k)}(\lambda)}
\newcommand\Tcl{T^{\flat}(\lambda)}
\newcommand\Fl{F(\lambda)}
\newcommand\BlVt{\tilde B_V(\lambda)}
\newcommand\KBlVt{K_{\BlVt}}
\newcommand\BlVaat{B_{V_{\alpha_1}}(\lambda)}
\newcommand\BV{B_V}
\newcommand\Bone{B_1}
\newcommand\Btwo{B_2}
\newcommand\Bthree{B_3}
\newcommand\Banyj{B_j}
\newcommand\PlV{P_V(\lambda)}
\newcommand\PlVc{P_V^{\flat}(\lambda)}
\newcommand\Pl{P_0(\lambda)}
\newcommand\SVl{S_V(\lambda)}
\newcommand\Sjr{S_j^{\reduced}}
\newcommand\Rkp{{\mathcal R}^k_+}
\newcommand\Rkm{{\mathcal R}^k_-}
\newcommand\Rkpm{{\mathcal R}^k_{\pm}}
\newcommand\Phys{{\mathcal P}}
\newcommand\Pc{\overline{\mathcal P}}
\newcommand\pip{\pi^{\perp}}
\newcommand\pipa{\pi_\partial}
\newcommand\gammapa{\gamma_\partial}
\newcommand\pipah{\hat\pi_\partial}
\newcommand\pit{\tilde\pi}
\newcommand\xit{\tilde\xi}
\newcommand\zetat{\tilde\zeta}
\newcommand\etat{\tilde\eta}
\newcommand\sigmat{\tilde\sigma}
\newcommand\sigmahat{\hat\sigma}
\newcommand\thetat{\tilde\theta}
\newcommand\psit{\tilde\psi}
\newcommand\phit{\tilde\phi}
\newcommand\chit{\tilde\chi}
\newcommand\rhot{\tilde\rho}
\newcommand\xib{\bar\xi}
\newcommand\zetab{\bar\zeta}
\newcommand\thetab{\bar\theta}
\newcommand\etab{\bar\eta}
\newcommand\iotal{\iota_{\lambda}}
\newcommand\rhoat{\rhot_{\alpha_1}}
\newcommand\Lambdat{\tilde\Lambda}
\newcommand\Lambdati{\tilde\Lambda^{\text{in}}}
\newcommand\Lambdato{\tilde\Lambda^{\text{out}}}
\newcommand\Lambdatp{\tilde\Lambda^{\text{prop}}}
\newcommand\Lambdai{\Lambda^{\text{in}}}
\newcommand\Lambdao{\Lambda^{\text{out}}}
\newcommand\poles{\Lambda'}
\newcommand\rpoles{\Lambda_p}
\newcommand\thresholds{\Lambda}
\newcommand\Vt{\tilde V}
\newcommand\It{\tilde I}
\newcommand\half{{\frac{1}{2}}}
\newcommand\sigmah{\sigma^{1/2}}
\newcommand\bX{\partial X}
\newcommand\bXb{\partial \Xb}
\newcommand\Deltabt{\tilde\Delta_0}
\newcommand\strip{\Omega_T}
\newcommand\Kf{K^{\flat}}
\newcommand\Gs{G^{\sharp}}
\newcommand\Gt{\tilde G}
\newcommand\Osc{\sci\Omega}
\newcommand\OSc{{}^\Scl\Omega}
\newcommand\Osch{\sci\Omega^{\half}}
\newcommand\Oscmh{\sci\Omega^{-\half}}
\newcommand\Isc{I_{sc}}
\newcommand\os{{\text{os}}}
\newcommand\Qzl{Q^0_{-\lambda}}
\newcommand\Lie{{\mathcal L}}
\newcommand\bl{{\text b}}
\newcommand\scl{{\text{sc}}}
\newcommand\sccl{{\text{scc}}}
\newcommand\Scl{{\text{sc}}}
\newcommand\ScLl{{\text{Sc,L}}}
\newcommand\ScRl{{\text{Sc,R}}}
\newcommand\Sccl{{\text{Scc}}}
\newcommand\sus{{\text{sus}}}
\newcommand\ssl{{\text{ss}}}
\newcommand\XXb{X^2_\bl}
\newcommand\XXbt{\Xt^2_\bl}
\newcommand\XXsc{X^2_\scl}
\newcommand\XXsct{\Xt^2_\scl}
\newcommand\XXSc{X^2_\Scl}
\newcommand\XXSct{\Xt^2_\Scl}
\newcommand\XXScL{X^2_\ScLl}
\newcommand\XXScR{X^2_\ScRl}
\newcommand\MMsc{M^2_\scl}
\newcommand\Deltab{\Delta_\bl}
\newcommand\Deltasc{\Delta_\scl}
\newcommand\DeltaSc{\Delta_\Scl}
\newcommand\DeltaScL{\Delta_\ScLl}
\newcommand\DeltaScR{\Delta_\ScRl}
\newcommand\prs{\sigma}
\newcommand\Nsc{N_\scl}
\newcommand\Nscp{N_{\scl,p}}
\newcommand\Nff{N_{\ff}}
\newcommand\Nffz{N_{\ff,0}}
\newcommand\Nffzp{N_{\ff,0,p}}
\newcommand\Nffl{N_{\ff,l}}
\newcommand\Nffml{N_{\ff,-l}}
\newcommand\Nmf{N_{\mf}}
\newcommand\Nmfz{N_{\mf,0}}
\newcommand\Nmfl{N_{\mf,l}}
\newcommand\Nmfml{N_{\mf,-l}}
\newcommand\ffb{\operatorname{bf}}
\newcommand\Ffb{\operatorname{bf'}}
\newcommand\ffsc{\operatorname{sf}}
\newcommand\ffSc{\operatorname{sf_C}}
\newcommand\Ffsc{\operatorname{sf'}}
\newcommand\rff{\rho_{\ff}}
\newcommand\rmf{\rho_{\mf}}
\newcommand\rffb{\rho_{\ffb}}
\newcommand\rffsc{\rho_{\ffsc}}
\newcommand\rFfsc{\rho_{\Ffsc}}
\newcommand\rffSc{\rho_{\ffSc}}
\newcommand\rinf{\rho_{\infty}}
\newcommand\CL{C_L}
\newcommand\CR{C_R}
\newcommand\betab{\beta_\bl}
\newcommand\betasc{\beta_\scl}
\newcommand\betaSc{\beta_\Scl}
\newcommand\BetaSc{\bar\beta_\Scl}
\newcommand\betaScL{\beta_\ScLl}
\newcommand\betaScR{\beta_\ScRl}
\newcommand\ScT{{}^\Scl T^*}
\newcommand\SccT{{}^\Scl \bar T^*}
\newcommand\ScS{{}^\Scl S^*}
\newcommand\Tb{{}^\bl T}
\newcommand\Tsc{{}^\scl T}
\newcommand\TSc{{}^\Scl T}
\newcommand\CSc{C_\Scl}
\newcommand\Lambdasc{{}^\scl\Lambda}
\newcommand\XXXb{X^3_\bl}
\newcommand\XXXsc{X^3_\scl}
\newcommand\XXXSc{X^3_\Scl}
\newcommand\XXXScO{X^3_{\Scl,O}}
\newcommand\XXXScF{X^3_{\Scl,F}}
\newcommand\XXXScS{X^3_{\Scl,S}}
\newcommand\XXXScC{X^3_{\Scl,C}}
\newcommand\KDsc{\operatorname{KD^{\half}_\scl}}
\newcommand\KDSc{\operatorname{KD^{\half}_\Scl}}
\newcommand\KDScEF{\operatorname{KD^{E,F}_\Scl}}
\newcommand\Oh{\operatorname{\Omega^{\half}}}
\newcommand\WFSc{\WF_\Scl}
\newcommand\WFtSc{\WF_{\text 3sc}}
\newcommand\WFScmf{\WF_{\Scl,\mf}}
\newcommand\WFScff{\WF_{\Scl,\ff}}
\newcommand\WFScs{\WF_{\Scl,\prs}}
\newcommand\WFScp{\WF'_\Scl}
\newcommand\WFScmfp{\WF'_{\Scl,\mf}}
\newcommand\WFScffp{\WF'_{\Scl,\ff}}
\newcommand\WFScsp{\WF'_{\Scl,\prs}}
\newcommand\Diffscc{\Diff_\sccl}
\newcommand\DiffSc{\Diff_\Scl}
\newcommand\DiffScc{\Diff_\Sccl}
\newcommand\DiffscI{\Diff_{\scl,\text{I}}}
\newcommand\VscI{\Vf_{\scl,\text{I}}}
\newcommand\DiffsV{\operatorname{Diff}_{\sus(V)}}
\newcommand\DiffsVsc{\operatorname{Diff}_{\sus(V),\scl}}
\newcommand\DiffsVCsc{\operatorname{Diff}_{\sus(V)-C,\scl}}   
\newcommand\Psisc{\Psop_\scl}
\newcommand\Psiscc{\Psop_\sccl}
\newcommand\Psiss{\Psop_\ssl}
\newcommand\Psisch{\Psop_{\scl,h}}
\newcommand\Psiscch{\Psop_{\sccl,h}}
\newcommand\PsiSc{\Psop_\Scl}
\newcommand\PsiScph{\Psop_{\Scl,\phi}}
\newcommand\PsiScra{\Psop_{\Scl,\rho^\sharp_a}}
\newcommand\PsiScc{\Psop_\Sccl}
\newcommand\PsiSccml{\Psop^{m,l}_\Sccl}
\newcommand\PsiScxx{\Psop^{*,*}_\Scl}
\newcommand\PsiScml{\Psop^{m,l}_\Scl}
\newcommand\PsiScmz{\Psop^{m,0}_\Scl}
\newcommand\PsiScmmz{\Psop^{-m,0}_\Scl}
\newcommand\PsiSckz{\Psop^{k,0}_\Scl}
\newcommand\PsiScmmml{\Psop^{-m,-l}_\Scl}
\newcommand\Psiscmkk{\Psop^{-k,k}_\scl}
\newcommand\Psiscmmmkk{\Psop^{-m-k,k}_\scl}
\newcommand\Psiscmoo{\Psop^{-1,1}_\scl}
\newcommand\Psiscmz{\Psop^{m,0}_\scl}
\newcommand\Psiscmmz{\Psop^{-m,0}_\scl}
\newcommand\PsiSckmkl{\Psop^{km,kl}_\Scl}
\newcommand\PsiScmplp{\Psop^{m',l'}_\Scl}
\newcommand\PsiScmmpllp{\Psop^{m+m',l+l'}_\Scl}
\newcommand\Psiscml{\Psop^{m,l}_\scl}
\newcommand\PsiScid{\Psop^{0,0}_\Scl}
\newcommand\PsiSczo{\Psop^{0,1}_\Scl}
\newcommand\PsiScmii{\Psop^{-\infty,\infty}_\Scl}
\newcommand\PsiScmiz{\Psop^{-\infty,0}_\Scl}
\newcommand\PsiScmoo{\Psop^{-1,1}_\Scl}
\newcommand\PsisCid{\Psop^{0,0}_{\scl-C}}
\newcommand\PsisC{\Psop_{\scl-C}}
\newcommand\Psiinf{\Psop_{\infty}}
\newcommand\Psiinfid{\Psop_{\infty}^0}
\newcommand\PsiFinf{\Psop_{\infty-\Fr}}
\newcommand\PsisVscml{\Psop^{m,l}_{\sus(V),\scl}}
\newcommand\PsisVsc{\Psop_{\sus(V),\scl}}
\newcommand\PsisVpsc{\Psop_{\sus(V_p),\scl}}
\newcommand\PsisVCSc{\Psop_{\sus(V)-C,\scl}}
\newcommand\SFinf{S_{\infty-\Fr}}
\newcommand\YsVC{Y^2_{\sus(V)-C,\scl}}
\newcommand\ffYsc{\ffsc_{\sus(V)}}
\newcommand\SXC{S(X;C)}
\newcommand\Ios{I_{\text{os}}}
\newcommand\pbL{\pi^2_{\bl,{\text L}}}
\newcommand\pbR{\pi^2_{\bl,{\text R}}}
\newcommand\pscL{\pi^2_{\scl,{\text L}}}
\newcommand\pscR{\pi^2_{\scl,{\text R}}}
\newcommand\PbO{\pi^3_{\bl,{\text O}}}
\newcommand\PscO{\pi^3_{\scl,{\text O}}}
\newcommand\PScO{\pi^3_{\Scl,{\text O}}}
\newcommand\PScF{\pi^3_{\Scl,{\text F}}}
\newcommand\PScC{\pi^3_{\Scl,{\text C}}}
\newcommand\PScS{\pi^3_{\Scl,{\text S}}}
\newcommand\pScL{\pi^2_{\Scl,{\text L}}}
\newcommand\pScR{\pi^2_{\Scl,{\text R}}}
\newcommand\CLF{\CL^F}
\newcommand\CLO{\CL^O}
\newcommand\CLS{\CL^S}
\newcommand\CLC{\CL^C}
\newcommand\DeltaYb{\Delta_{\bl,Y}}
\newcommand\DeltaYsc{\Delta_{\sus-\scl}}
\newcommand\diag{\operatorname{diag}}
\newcommand\Vf{{\mathcal V}}
\newcommand\Vb{{\mathcal V}_{\bl}}
\newcommand\Vsc{{\mathcal V}_{\scl}}
\newcommand\VSc{{\mathcal V}_{\Scl}}
\newcommand\VfI{\Vf_{\text{I}}}
\newcommand\VfIq{\Vf_{\text{I},q}}
\newcommand\scH{{}^\scl H}
\newcommand\scHg{\scH_g}
\newcommand\Hss{H_\ssl}
\newcommand\xh{\hat x}
\newcommand\yh{\hat y}
\newcommand\sh{\hat s}
\newcommand\rh{\hat r}
\newcommand\Yh{\hat Y}
\newcommand\Zh{\hat Z}
\newcommand\Yb{\bar Y}
\newcommand\hb{\bar h}
\newcommand\xih{\hat\xi}
\newcommand\etah{\hat\eta}
\newcommand\muh{\hat\mu}
\newcommand\mub{\bar\mu}
\newcommand\nub{\bar\nu}
\newcommand\mubh{\widehat{\bar\mu}}
\newcommand\yb{\bar y}
\newcommand\zb{\bar z}
\newcommand\ub{\bar u}
\newcommand\Qb{\bar Q}
\newcommand\Wbp{{\bar W}^\perp}
\newcommand\Wp{W^\perp}
\newcommand\Kt{\tilde K}
\newcommand\Wt{\tilde W}
\newcommand\Ut{\tilde U}
\newcommand\yt{\tilde y}
\newcommand\ut{\tilde u}
\newcommand\vt{\tilde v}
\newcommand\ft{\tilde f}
\newcommand\htil{\tilde h}
\newcommand\St{\tilde S}
\newcommand\Pt{\tilde P}
\newcommand\Rt{\tilde R}
\newcommand\qt{\tilde q}
\newcommand\Qt{\tilde Q}
\newcommand\Xb{\bar X}
\newcommand\lambdat{\tilde\lambda}
\newcommand\betat{\tilde\beta}
\newcommand\Phit{\tilde\Phi}
\newcommand\epst{\tilde\epsilon}
\newcommand\ep{\epsilon}
\newcommand\bt{\tilde b}
\newcommand\Xt{\tilde X}
\newcommand\Mt{\tilde M}
\newcommand\At{\tilde A}
\newcommand\Et{\tilde E}
\newcommand\Ht{\tilde H}
\newcommand\at{\tilde a}
\newcommand\Ct{\tilde C}
\newcommand\pih{\hat\pi}
\newcommand\Rh{\hat R}
\newcommand\Ah{\hat A}
\newcommand\Bh{\hat B}
\newcommand\Ch{\hat C}
\newcommand\Gh{\hat G}
\newcommand\Hh{\hat H}
\newcommand\Qh{\hat Q}
\newcommand\Ph{\hat P}
\newcommand\Nh{\hat N}
\newcommand\Sh{\hat S}
\newcommand\Gcal{{\mathcal G}}
\newcommand\GcalC{{\mathcal G}_C}
\newcommand\Jcal{{\mathcal J}}
\newcommand\JcalC{{\mathcal J}_C}
\setcounter{secnumdepth}{3}
\newtheorem{lemma}{Lemma}[section]
\newtheorem{prop}[lemma]{Proposition}
\newtheorem{thm}[lemma]{Theorem}
\newtheorem{cor}[lemma]{Corollary}
\newtheorem{result}[lemma]{Result}
\newtheorem*{thm*}{Theorem}
\newtheorem*{prop*}{Proposition}
\newtheorem*{cor*}{Corollary}
\newtheorem*{conj*}{Conjecture}
\numberwithin{equation}{section}
\theoremstyle{remark}
\newtheorem{rem}[lemma]{Remark}
\newtheorem*{rem*}{Remark}
\theoremstyle{definition}
\newtheorem{Def}[lemma]{Definition}
\newtheorem*{Def*}{Definition}
\def\signature#1#2{\par\noindent#1\dotfill\null\\*
{\raggedleft #2\par}}

\renewcommand{\theenumi}{\roman{enumi}}
\renewcommand{\labelenumi}{(\theenumi)}

\title{Geometry and analysis in many-body scattering}
\author{Andr\'as Vasy}
\address{Department of Mathematics, Massachusetts Institute of Technology,
Cambridge MA 02139, U.S.A.}
\email{andras@math.mit.edu}
\thanks{This work is partially supported by NSF grant \#DMS-0201092,
a Fellowship from the Alfred P.\ Sloan Foundation, and
the Universit\'e de Nantes, where
these lectures were originally given.}

\maketitle

\section{Introduction}
The present notes are an effort to explain in relatively non-technical
terms recent results in many-body scattering and related topics.
Thus, many results in the many-body setting should be understood
as new results on the propagation of singularities, here understood
as lack of decay of wave functions
at infinity, with much in common with real principal
type propagation, i.e.\ wave phenomena. Motivated by this, I first
briefly describe propagation of singularities for the wave equation.
This is a remarkable relationship between geometric optics (the
particle view of light) and the solutions of the wave equation
(the wave view).

Next, in Section~\ref{sec:mb},
I explain the geometry of many-body scattering, which includes
both that of the configuration space and phase space. This geometry
is closely related to classical mechanics, playing the role
of geometric optics, but even at this point quantum phenomena
emerge. This leads to the analytic results, namely the propagation
of singularities connecting classical and quantum mechanics.

Much as for the wave equation, such a result has immediate applications,
including the description of the scattering matrices and of
the scattering phase. Slightly stronger
versions can even lead to inverse results, a topic covered in the
following section.

After so explaining the results, in Sections~\ref{sec:psdo}-\ref{sec:positive},
I will try to at least give a flavor
of how they are proved. This uses a many-body pseudo-differential
algebra and positive commutator estimates, so these are discussed.
We remark that these techniques are closely related to the proofs
of the propagation of singularities for the wave equation, but there
are significant differences as well, mostly arising from bound states
of particles, which have no analogues for the wave equation.
The pseudo-differential algebra
itself is very interesting from the viewpoint of non-commutative
geometry: there is a hierarchy of operator valued symbols at infinity.

Asymptotic completeness was the main focus of work in
many-body scattering for a long period. In Section~\ref{sec:ac}, I briefly
explain how it relates to the microlocal estimates.

There is another area that is very closely related to many-body scattering,
namely scattering on higher rank non-compact symmetric spaces. Here,
in Section~\ref{sec:sl3}, we only
discuss rank two, which corresponds to three-body scattering,
since this is the only part that has been properly
written up, but it is expected that very soon these results will extend
to all higher rank spaces.

I hope that the notes will make many of these results more accessible,
the connections more transparent, and
explain the motivation behind them. Many-body scattering has a long
history, and here I can only talk about the most recent developments.
An excellent overview of results known in the early 1990s can be
found in Hiroshi Isozaki's lecture notes \cite{Iso}. Indeed, in
some sense, the current notes continue where \cite{Iso} left off.
I introduce a fully microlocal picture, motivated by the geometric
approach of Richard Melrose \cite{RBMGeo}, and emphasize the results
these give, but the basic spectral and scattering results follow
from a simpler `partial' microlocalization, which is one of
the subjects of \cite{Iso}.

The notes were originally prepared for a mini-course at the
Universit\'e de Nantes at the invitation of Professor Xue-Ping Wang,
whose hospitality I gratefully acknowledge. The analytic continuation
of the resolvent
on symmetric spaces is a more recent development, but it was fueled
by a discussion during the visit of Rafe Mazzeo, my collaborator, to Nantes.
Over an espresso, Gilles Carr\'on mentioned that the existence of the analytic
continuation was not known, something that was hard to believe, but we
immediately realized that our methods should yield such a continuation
rather directly. I also thank Gunther Uhlmann for urging me to write up
these notes: without him, they may never have been written up, and Rafe
Mazzeo for a careful reading of the manuscript.

\section{Geometric optics and the wave equation}
According to the rules of geometric optics, light propagates in straight
lines, and reflects/refracts from surfaces according to Snell's/Descartes' law.
That is to say, considering light as a stream of billiard balls, the
energy as well as the tangential component of the momentum (tangential to
the surface hit) is conserved upon hitting the surface.

But light satisfies the wave equation, i.e.\ if $u=u(x,t)$
is the electromagnetic
field on $\Omega_x\times\Real_t$, $\Omega\subset\Real^n$,
then $Pu=0$ where $P$ is the wave operator
$c^2\Delta-D_t^2$, and a boundary condition also holds (say, Dirichlet),
if $\Omega$ is not the whole space. (Here $D_t=\frac{1}{i}\,\frac{\pa}{\pa t}$
and $\Delta=\sum_j D_{x_j}^2$ is the {\em positive} Laplacian.)
How are these two viewpoints related?

One can phrase the connection in different ways. The most usual one in physics
is that the billiard ball picture is accurate in the high frequency, i.e.\ low
wave length, limit. That is to say, for high frequency light, geometric optics
is accurate up to a `small' error. A slightly different way of looking at
this, which however does not involve approximations, is that the location
of singularities of the solution of the wave equation is {\em exactly}
predicted by geometric optics. Here singularities are understood as
lack of smoothness, or possibly lack of analyticity.

Indeed, it is convenient at this point to generalize the setting somewhat.
So let $(\Omega,g)$ be a Riemannian manifold with corners, $P=c^2\Delta_g
-D_t^2$, $c>0$. The speed of light, $c$,
may be absorbed in the metric $g$, of course,
we keep the notation in analogy with the usual wave equation.

For simplicity of notation in this paragraph
we assume that $M_z=\Omega_x\times\Real_t$ is boundaryless; in general,
the same definitions hold in the interior of $M$.
Thus, we associate a homogeneous real function on $T^*M$
to $P$, namely its principal symbol: $p=c^2|\xi|_g^2-\tau^2$,
where we write $\zeta=(\xi,\tau)$ as the dual variable of $z=(x,t)$.
Now $T^*M$ is a symplectic manifold with symplectic form $\omega=\sum
d\zeta_j\wedge dz_j$. Thus, $p$ gives rise to a vector field $H_p$, called
the Hamilton vector field, by requiring that
$\omega(V,H_p)=Vp$ for any vector field
$V$ on $T^*M$. Hence $H_p$ is a smooth vector field on $T^*M$ explicitly
given by
\begin{equation*}
H_p=\frac{\pa p}{\pa\zeta}\,\frac{\pa}{\pa z}-\frac{\pa p}{\pa z}\,\frac{\pa}
{\pa\zeta}.
\end{equation*}
Note that $p$ is constant along the integral curves of $H_p$ since
taking $V=H_p$, $0=\omega(H_p,H_p)=H_p p$. Null bicharacteristics are
the integral curves of $H_p$ inside its characteristic set
$\Sigma=p^{-1}(\{0\})$. Thus, if $\gamma:I\to T^*M$ is a null bicharacteristic
(here $I$ is an interval), and $z(s)=z(\gamma(s))$,
$\zeta(s)=\zeta(\gamma(s))$, then these solve the ODE's
$\frac{dz}{ds}=\frac{\pa p}{\pa\zeta}$, $\frac{d\zeta}{ds}=
-\frac{\pa p}{\pa z}$. Hence, when $M=\Omega\times\Real$,
$\Omega\subset\Real^n$, $p=c^2|\xi|^2-\tau^2$ as above, we deduce that $\xi$
and $\tau$ are constant along the integral curves of $H_p$, hence their
projection to $M$ consists of straight line segments. More generally,
the projection of null-bicharacteristics to $\Omega$ are geodesics of $g$.

There is an appropriate extension of this at boundary surfaces and even at
corners, called generalized broken bicharacteristics, see
\cite{Melrose-Sjostrand:I, Lebeau:Propagation},
which I will not explain in full generality,
though I remark that many-body scattering, discussed in the next section
in detail,
is rather similar. However, a somewhat typical example is that of broken
bicharacteristics. These are piecewise bicharacteristics, i.e.\ there
is a sequence $s_j$, $j$ in a subset of integers, such that for each $j$,
$\gamma|_{(s_j,s_{j+1})}$ is a null bicharacteristic in the sense
described above, the projection $z\circ\gamma$
of $\gamma$ to $M$ is continuous,
and $\gamma(s_j+)-\gamma(s_j-)$ is conormal
to the smallest dimensional
boundary face containing $z(\gamma(s_j))$. Thus, the tangent
vectors to $z\circ\gamma|_{(s_j,s_{j+1})}$ and
$z\circ\gamma|_{(s_{j-1},s_j)}$ differ by a vector normal to
the smallest boundary face containing $z(\gamma(s_j))$. This expresses that the
normal component of the momentum may change, while the tangential component
is conserved, when a light ray hits a boundary.

\begin{figure}
\begin{center}
\mbox{\epsfig{file=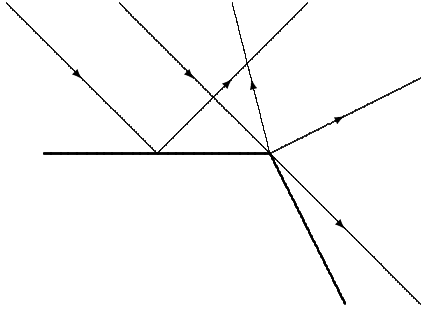}}
\end{center}
\caption{Projection of broken bicharacteristics to $\Omega$. When
rays hit the boundary hypersurfaces, the tangential component of the
momentum and the kinetic energy are conserved, but the normal
component may change. At the corner, there is no tangential component
(though there would be if the time variable were not projected out),
so the only constraint is the conservation of kinetic energy.}
\label{fig:brray}
\end{figure}

Now one can describe the singularities of $u$ using null bicharacteristics.
Let $o$ be the zero section of $T^*M$.
The location of the singularities is described by an object
\begin{equation*}
\WF(u)\subset T^*M\setminus o=\{(z,\zeta):\ \zeta\neq 0\}
\end{equation*}
that is conic in $\zeta$, i.e.\ $(z,\zeta)\in\WF(u)$
if and only if $(z,r\zeta)\in\WF(u)$ for every $r>0$. $\WF(u)$ is
called the wave front set of $u$, and it describes where (in $z$)
and in which codirection $\zeta$ is the distribution $u$ not $\Cinf$.
More precisely, the definition of $\WF(u)$ is that $(z_0,\zeta_0)\nin
\WF(u)$ if and only if there exists $\phi\in\Cinf_c(M)$, $\phi(z_0)\neq 0$
such that the Fourier transform $\Fr(\phi u)$ of $\phi u$
is rapidly decreasing in an open cone around $\zeta_0$. Here we assume
that $M$ is boundaryless; otherwise we need to require that $\phi$ is
supported in the interior of $M$. Again, there is a natural definition at
$\pa M$ which we do not give here. (There are
more natural versions of this definition using pseudo-differential
operators that I will describe later.) As an example, consider the
step function: write $z=(z_1,z'')$, $u(z)=1$ if $z_1>0$, $u(z)=0$ if
$z_1<0$. Then
\begin{equation*}
\WF(u)=N^*\{z_1=0\}\setminus o=\{(0,z'',\zeta_1,0):\ \zeta_1\neq 0\},
\end{equation*}
the conormal bundle of the hypersurface $z_1=0$, with its zero section
removed. The same statement holds, with $=$ possibly replaced by
$\subset$, if we take any $\Cinf$ function $u_0$ on $M$,
and then define $u=u_0$ in $z_1>0$ and $u=0$ in $z_1<0$.
Informally, one might say that $u$ is singular in $z_1$ at $z_1=0$,
but it depends smoothly on $z''$. The wave front set thus pinpoints
not only the locations $z$ of singularities (lack of smoothness) in $M$,
but it refines it by also giving the frequencies (or rather direction of
frequencies) at which these appear at $z$.

The theorem we are after is the following. In early versions it goes
back to Lax \cite{Lax:Asymptotic}, its boundaryless version is due to
H\"ormander \cite{Hormander:Existence},
the smooth boundary versions are due to Melrose, Sj\"ostrand, Taylor and Ivrii
\cite{Ivrii:Wave, Melrose-Sjostrand:I, Melrose-Sjostrand:II, Taylor:Grazing},
and the corner version in the analytic category is due to
Lebeau~\cite{Lebeau:Propagation} (the
$\Cinf$ version is still not known in the corner setting) while a different
extension, to conic
points, is due to Melrose and Wunsch \cite{Melrose-Wunsch:Propagation}.

\begin{thm}
Suppose $Pu\in\Cinf(M)$, and if $\pa M\neq 0$ then $u|_{\pa M}=0$.
Then $\WF(u)\subset\Sigma=p^{-1}(\{0\})$
(microlocal elliptic regularity). Moreover, $\WF(u)$ is a union of maximally
extended generalized broken bicharacteristics inside $\Sigma$ (propagation
of singularities).
\end{thm}

This theorem states that if a point $(z,\zeta)\in T^*M\setminus o$
is in $\WF(u)$ and $u$ solves $Pu\in\Cinf(M)$, and satisfies a
boundary condition if appropriate, then there is at least one
maximally extended generalized broken bicharacteristic through $(z,\zeta)$
that is completely contained in $\WF(u)$. Of course, in the absence
of boundaries, and often even in their presence, there is a unique
maximally extended generalized broken bicharacteristic through $(z,\zeta)$,
so the statement is that this bicharacteristic is completely in $\WF(u)$.
However, as soon as codimension two or higher corners appear, there
is no hope for such uniqueness, and this theorem is the optimal
statement.

At least in the nicest settings (no boundaries, or non-degeneracy assumption
at the boundaries which are assumed to be smooth), this theorem can
be improved significantly to predict not only the location, but also
the amplitude of the singularities of $u$.

\section{Propagation in many-body scattering}\label{sec:mb}
There is an analogous setup for scattering. Now we want to understand how
interacting particles behave. Again, there is a classical mechanical setup
(the analogue of geometric optics) and a quantum mechanical setup
(the analogue of the wave equation). To focus on the most relevant points,
I formulate the problem in a time-independent fashion, though it is
easy to reformulate everything in a time dependent way. We only
do this in a remark following Theorem~\ref{thm:prop-sing}.

Thus, we want to understand tempered distributional solutions $u$ of
$(H-\lambda)u=0$; here $\lambda\in\Real$ is the energy, and $H$ is the
Hamiltonian, i.e.\ the analogue of $H-\lambda$ is $P$ above. Namely, if we have
$N$ particles, each of which is $d$-dimensional with positions $x_1,
\ldots, x_N\in\Real^d$, mass $m_1,\ldots,m_N$, and the interaction
between particle $i$ and $j$ is given by a potential $V_{ij}$ (which
is a function on $\Real^d$), then the Hamiltonian describing this
system is
\begin{equation*}
H=\sum_{i=1}^N\frac{1}{2m_i}\Delta_{x_i}+\sum_{i<j} V_{ij}(x_i-x_j)=\Delta+V,
\end{equation*}
which is an operator on (functions on) $\Real^n=\Real^{Nd}$. Planck's
constant $\hbar$ is here taken to be $1$; it could be absorbed in the $x_i$
by a simple rescaling.

Now $H$ is elliptic in the standard sense, namely its principal symbol is
$\sum\frac{1}{2m_i}|\xi|^2$, which never vanishes outside the zero section $o$.
Note that the potential is lower order than $\Delta$ in the standard sense,
so it is not part of the principal symbol.
So, by the previous theorem,
\begin{equation*}
(H-\lambda)u=0\Rightarrow\WF(u)=\emptyset\Rightarrow u\in\Cinf(\Real^n).
\end{equation*}
So the only possibility of interesting behavior for $u$ is at infinity,
and this is exactly what we want to understand.

\begin{figure}
\begin{center}
\mbox{\epsfig{file=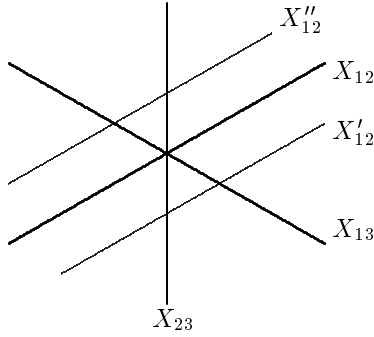}}
\end{center}
\caption{Collision planes $X_{12}$, $X_{13}$ and $X_{23}$ and translates
$X'_{12}$ and $X''_{12}$ of $X_{12}$. $V_{12}$ is constant along $X_{12}$, is a
(typically different) constant along $X'_{12}$, etc., so it does not
decay at infinity unless it is identically zero.}
\label{fig:collplij}
\end{figure}

The main feature of many-body problems is that even if $V_{ij}$
decays at infinity
on $\Real^d$, it does {\em not} decay at infinity in $\Rn$ since it is
a constant along $X_{ij}=\{x_i=x_j\}$,
as well as along its translates $X_{ij}'$, $X_{ij}''$,
so it does not decay if we go to infinity, say, along $X_{ij}$; see
Figure~\ref{fig:collplij}.
The $X_{ij}$ are called collision planes (as are their intersections)
since at $X_{ij}$ particles $i$ and $j$ are at the same place.

In the two-body problem one actually has $H=\Delta_{x_1,x_2}+V_{12}(x_1-x_2)$,
i.e.\ $V_{12}$ still does not decay at infinity, e.g.\ if one keeps
$x_1=x_2$ but lets $x_1\to\infty$. However, one can easily remove the
center of mass by performing a Fourier transform along $X_{12}$.
This conjugates $H-\lambda$ to $H^{12}+|\xi_{12}|^2-\lambda$, where
$\xi_{12}$ is the variable on $X_{12}^*$, and $H^{12}=\Delta_{X^{12}}+V_{12}$,
$X^{12}$ being the orthocomplement of $X_{12}$.
Thus, one reduces the study of $H-\lambda$ to that
of a Hamiltonian on $X^{12}$, namely $H^{12}-\lambda'$,
$\lambda'=\lambda-|\xi_{12}|^2$ being a shifted spectral parameter.
Now $V_{12}$ decays at infinity (we are working
on $X^{12}$!), so $H^{12}$ can be considered as
a perturbation of $\Delta_{X^{12}}$, hence its analysis is rather simple.
Notice that the point spectrum of $H^{12}$ gives rise to a branch of
the continuous spectrum of $H$: this is a phenomenon that is very
typical in many-body scattering.
The center of mass can also be removed in any actual many-body problem,
but one still obtains a Hamiltonian with non-decaying potentials as before.

One can still talk about classical mechanics, just as for the wave equation,
using bicharacteristics. These are deterministic -- if $V$ is smooth enough
(we usually assume that $V$ is $\Cinf$). But much like for corners, there
is a compressed description of dynamics near infinity. This is
somewhat more complicated than for the wave equation, but only
because particles can be bound together. Thus, even the `classical'
description is partly quantum. These two facts, the presence of collision
planes and the bound states, are the two crucial features of many-body
scattering.

The compressed dynamics in the absence of bound states looks just like
in the wave equation setting. One should think of this as a good
description when a classical trajectory is uniformly near infinity.

\begin{figure}
\begin{center}
\mbox{\epsfig{file=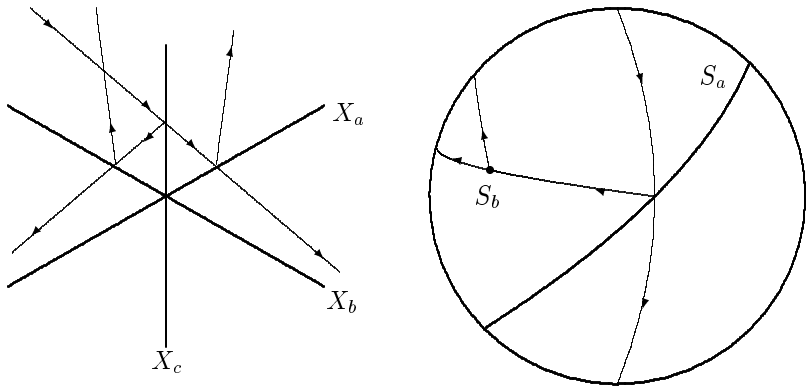}}
\end{center}
\caption{On the left,
broken geodesics in $\Real^n\setminus\{0\}$, $n=2$,
broken at the collision planes
$X_a$, $X_b$ and $X_c$. On the right, the projection of broken
geodesics in $\Real^n\setminus\{0\}$, $n=3$, emanating from the north pole,
to the unit sphere $S_0$, better understood
as the sphere at infinity. The $C_a$, $C_b$ are the intersection
of the collision planes $X_a$, $X_b$ with $S_0$; $\dim X_a=2$, $\dim X_b=1$.}
\label{fig:brtraj}
\end{figure}

More precisely, it is convenient to introduce Agmon's generalization of
the many-body problem, which amounts to using the vector space structure
of $\Real^n$ as the setting. One can also give geometric generalizations
(in the sense of differential geometry)
that arose from the work of Melrose \cite{RBMSpec, RBMGeo},
and I will do this later.

So we work on the vector space $X_0=\Rn$, equipped with the Euclidean
metric. We are also given a finite collection $\calX=\{X_a:\ a\in I\}$
of linear subspaces $X_a$ of $\Rn$, called the collision planes.
We assume that $\calX$ is closed under intersections, and
$X_0=\Rn\in\calX$, $X_1=\{0\}\in\calX$. We
let $X^a=X_a^\perp$ be the orthocomplement of $X_a$ in $\Rn$, so
$\Rn=X_a\oplus X^a$.
(Agmon's generalization is thus that the $X_a$ do not have to come from
intersections of the planes $X_{ij}=\{x_i=x_j\}$.)
We write the corresponding coordinates as $(x_a,x^a)$,
and denote the orthogonal projection to $X^a$ by $\pi^a$. A many-body
Hamiltonian in potential scattering is an operator of the form
\begin{equation*}
H=\Delta+\sum_a (\pi^a)^*V_a,
\end{equation*}
where $V_a$ is a real valued function on $X^a$ in a certain class,
for example $V_a$ is a symbol on $X^a$ of negative order:
$V_a\in S^{-\rho}(X^a)$, $\rho>0$. We also assume that $V_0=0$
for normalization; note that $X^0=\{0\}$,
so $V_0$ would simply play the role of the spectral parameter.
We sometimes drop the pull-back notation
from now on and write $H=\Delta+\sum_a V_a$.

Another useful piece of terminology is the following. We say that $V_a$ is
short range if $V_a\in S^{-\rho}(X^a)$ for some $\rho>1$. We say
that $V_a$ is long-range if $V_a\in S^{-\rho}(X^a)$ for some
$\rho\in(0,1]$. The Coulomb potential is thus `marginally long-range',
at least if we ignore its singularity at $0$ (which is not a serious
problem anyway). Whether $V_a$ is short- or long-range does not make
any difference for the propagation phenomena we discuss in this section.
However, it does make a major difference for the precise behavior of
generalized eigenfunctions at the `radial sets' which we discuss later.
This also shows up in the related issue of asymptotic completeness.

Yet another notation we use on occasion is that of a $k$-cluster.
Physically, a cluster describes particles that are close (or collide),
and a $k$-cluster
means that there are $k$ clusters of particles, inside each of which
the particles
are close to each other. So in $N$-particle scattering, the $N$-cluster
describes $N$ asymptotically free particles (none is close
to any other), hence we say that the
collision plane $X_0=\Rn$ is the $N$-cluster. On the other hand,
if $X_a\neq \{0\}$ is such that $X_b\subsetneq X_a$ implies that $X_b=\{0\}$,
then $X_a$, or rather $a$, is a 2-cluster. E.g.\ given five particles,
a 2-cluster is where $x_1=x_2$ and $x_3=x_4=x_5$, i.e. the particles
1 and 2, resp. 3, 4 and 5, are close to each other. In general, a $k$-cluster
$X_a$ can be defined by the length of nested chains of collision planes inside
$X_a$.

One need not assume that all interactions between the particles are via
potentials. Indeed, $V_a$ may be allowed to be any first order
differential operator on the vector space $X^a$ with symbolic coefficients
of negative order. Also, one may generalize the metric $g$ in an analogous
fashion, as discussed later, which in effect allows $V_a$ to be second
order provided that $H$ remains elliptic. To simplify the notation,
and due to the traditions, we mostly talk as if $V_a$ were potentials,
but the generalization to such higher order perturbations requires only
occasional and minor modifications, which will be pointed out.

The subsystem Hamiltonians are defined by
\begin{equation*}
H^a=\Delta_{X^a}+\sum_{X_a\subset X_b} V_b.
\end{equation*}
Note that $X_a\subset X_b$ if and only if $X^a\supset X^b$, so above $V_b$
is really the pull-back of $V_b$ from $X^b$ to to $X^a$ by the
orthogonal projection. Thus, $H^a$ is an operator on (functions on)
$X^a$, and indeed it is a many-body Hamiltonian.

We also let
\begin{equation*}
X_{a,\sing}=\bigcup\{X_b:X_b\subsetneq X_a\},
\ X_{a,\reg}=X_a\setminus X_{a,\sing},
\end{equation*}
be the singular, respectively the regular, part of $X_a$.
Thus, if $X_c$ is a collision plane and $X_a$ is not a subset of $X_c$,
then $X_a\cap X_c$ is a proper subset of $X_a$, and is a collision plane
(since $\calX$ is closed under intersections), so $X_a\cap X_c\subset
X_{a,\sing}$. Correspondingly, $V_c$ decays at $X_{a,\reg}$, so
\begin{equation*}
H_a=\Delta_{X_a}+H^a,
\end{equation*}
which is an operator on (functions on) $X_0=\Rn$, has the property that
$H-H_a$ is a function that decreases at $X_{a,\reg}$. So $H_a$
should be thought of as a good approximation of $H$ at $X_{a,\reg}$.
Note that $X_{a,\sing}$ is a finite union of codimension $\geq 1$ submanifolds
of $X_a$, so $X_{a,\reg}$ is in particular an open dense subset of $X_a$.
Also, note that $\Delta_{X_a}$ plays a role analogous to the kinetic
energy of the center of mass in the two-body setting, but now this description
only valid locally, at $X_{a,\reg}$.

Having thus described the configuration space $X=X_0=\Rn$, the next step is
to describe the phase space, as was done first in \cite{Vasy:Propagation-Many}
and \cite{Vasy:Bound-States}. The main goal in the process is to obtain
a space on which broken bicharacteristics behave well. We remind the
reader that we are concerned with singularities at infinity, hence
with bicharacteristics that are uniformly close to infinity. Later
we give a compactified description, but here for simplicity we give its
homogeneous version, much as for the wave equation where bicharacteristics were
integral curves of the homogeneous principal symbol. So we start with $T^*X$,
but we wish to compress it at $X_a$ in such a way that at $X_{a,\reg}$,
$T^*_{X_{a,\reg}}X$ is replaced by $T^*_{X_{a,\reg}}X_a=T^*X_{a,\reg}$.
For broken
bicharacteristics this has the effect that only the $X_a$-tangential component
of the momentum is preserved at $X_{a,\reg}$. So we define the
compressed cotangent bundle as
\begin{equation*}
\dT X=\cup_{a\neq 1}T^*X_{a,\reg}.
\end{equation*}
Note that this is at first just a set, equipped with a projection $\dT X\to
X\setminus\{0\}$ induced by the bundle projections $T^*X_{a,\reg}
\to X_{a,\reg}$. There is also a natural $\Real^+$-action on $\dT X$ via
dilation in the configuration variables:
\begin{equation}\label{eq:config-dilation}
\Real_r^+\times T^*X_{a,\reg}\ni (r,x_a,\xi_a)\mapsto
(rx_a,\xi_a)\in T^*X_{a,\reg}.
\end{equation}
We topologize $\dT X$ via the projection
\begin{equation*}
\pi: T^*_{X\setminus\{0\}}X\to\dT X,
\end{equation*}
whose restriction to $T^*_{X_{a,\reg}}X$ is the pull-back of one-forms
by the inclusion map $X_{a,\reg}\hookrightarrow X$. Thus, writing
$(\xi_a,\xi^a)$ as the momenta dual to $(x_a,x^a)$, $\pi$
projects out the normal component of the momentum, $\xi^a$. The topology
is then the weakest topology that makes $\pi$ continuous, i.e.\ a set $C$
in $\dT X$ is closed if and only if $\pi^{-1}(C)$ is closed.

We can now describe the contribution of the bound states to the characteristic
sets. As mentioned above, this is one of the most interesting features
of many-body scattering that has no analogue for the wave equation.
These are conic subsets of $\dT X$ (conic with respect to the
$\Real_r^+$-action in \eqref{eq:config-dilation}). The characteristic
sets describe where certain operators are not elliptic, i.e.\ invertible,
at infinity, in a precise sense described in the subsequent sections.
They correspond to the `energy shell', i.e.\ being on the characteristic
set at energy $\lambda$ means that the particles have total energy $\lambda$.
We let
\begin{equation*}
\Char_0(\lambda)=\{(x,\xi)\in T^*X:\ g(\xi)=\lambda\}
\end{equation*}
be the free characteristic variety, with $g$ being the metric function
on $T^*X$, and more generally we set
\begin{equation*}
\Char_a(\lambda)=\{(x_a,\xi_a)\in T^*X_a:\ \lambda-g_a(\xi_a)\in
\pspec(H^a)\}\subset T^*X_a.
\end{equation*}
Notice that $\lambda=g_a(\xi_a)+\ep_\alpha$,
$\ep_\alpha\in \pspec(H^a)$, corresponds to the splitting of the total energy
$\lambda$ to the kinetic energy of the cluster, $g_a(\xi_a)$, plus the
energy of the bound state, $\ep_\alpha$. Thus, $\Char_a(\lambda)$
describes that particles may exist in a bound state of $H^a$,
of energy $\ep_\alpha$, along $X_a$, with kinetic energy $g_a(\xi_a)=\lambda-
\ep_\alpha$.
Moreover, $H^0$ is the zero operator on $X^0=\{0\}$, so if $a=0$, these
two definitions are consistent.
If $X_a\subset X_b$, the pull-back of one-forms gives a projection
$\pi_{ba}:T^*_{X_a}X_b\to T^*X_a$. Let
\begin{equation*}\begin{split}
&\dChar(\lambda)=\cup\dChar_a(\lambda)\subset\dT X,\\
&\dChar_a(\lambda)=\cup_{X_b\supset X_a} \pi_{ba}(\Char_b(\lambda))
\cap T^*X_{a,\reg}.
\end{split}\end{equation*}

\begin{figure}
\begin{center}
\mbox{\epsfig{file=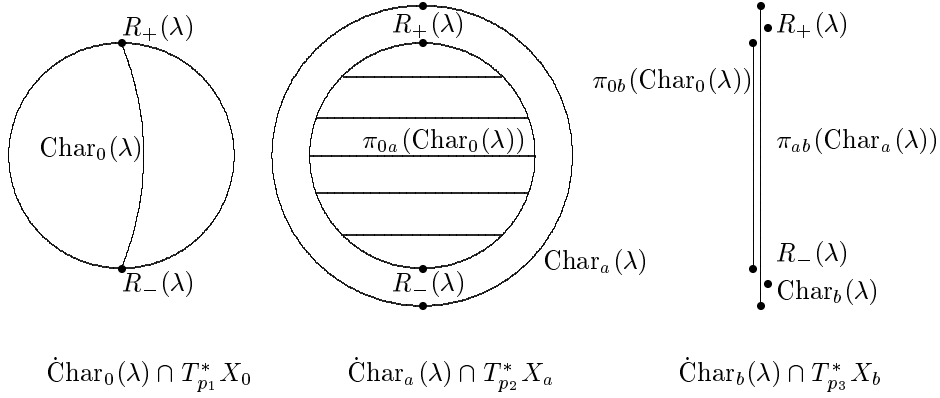}}
\end{center}
\caption{The characteristic set of many-body Hamiltonians. Here
$H$ is a 4-body Hamiltonian, $a$ is a 3-cluster, $b$ is a 2-cluster,
$p_1\in X_{0,\reg}$, $p_2\in X_{a,\reg}$, $p_3\in X_{b,\reg}$.
The solid dots are the radial sets, defined below.}
\label{fig:charset}
\end{figure}

In order to understand $\dChar(\lambda)$ it is important to keep in mind
several results on the structure of the eigenvalues of the subsystems. So let
\begin{equation*}
\Lambda_a=\cup_{b:X^b\subsetneq X^a}\pspec(H^b)
\end{equation*}
be the set of thresholds of $H^a$. Fundamental results of
Perry, Sigal and Simon \cite{Perry-Sigal-Simon:Spectral} and of
Froese and Herbst \cite{FroMourre} show that $\Lambda_a$ is closed,
countable, and the countable set $\pspec(H^a)$ can only accummulate
at $\Lambda_a$, so
\begin{equation*}
\Lambda'_a=\Lambda_a\cup\pspec(H^a)=\cup_{b:X^b\subset X^a}\pspec(H^b)
\end{equation*}
is also closed. Hence, $\dChar(\lambda)$ is a closed subset of $\dT X$.
In fact, the quotient $\dot\Sigma(\lambda)$ of
$\dChar(\lambda)$ by the $\Real^+$ action (which can be realized by
restricting the various bundles
to the unit sphere, $S_0=\{x\in X_0:\ |x|=1\}$) is compact, and indeed
it is metrizable, see \cite{Vasy:Propagation-Many}. Since compact
topological spaces have better properties than non-compact ones, it
is quite natural to work with $\dot\Sigma(\lambda)$, although we
do not follow this route in this section. We also remark that it is much
better to think of $\dot\Sigma(\lambda)$ lying at the sphere at infinity,
rather than at $S_0$, since it is the dynamics at infinity that is described
here. We will take up this approach in later sections.

We also recall another result of Froese and Herbst \cite{FroExp}, namely
that eigenfunctions $\psi_\alpha$ of $H^a$ with eigenvalue $\ep_\alpha$
decay exponentially on $X^a$, at a specified rate, if $\ep_\alpha\nin
\Lambda_a$. This generalizes to higher order perturbations, but requires
a somewhat different approach, see \cite{Vasy:Exponential}. In fact, this
is the only place where second order perturbations behave differently
from first or zeroth order ones. For the latter, there can be no
positive energy bound states, while for the former this has been only
proved for small metric perturbations in \cite{Vasy:Exponential},
and it is not clear whether it
holds more generally, especially for trapping perturbations. (Note that if
$H^a$, $a\neq 1$, is trapping then $H$ is trapping at infinity!)

A generalized broken bicharacteristic (at energy $\lambda$)
is then a continuous map $\gamma:I\to\dChar(\lambda)$, $I$ an interval,
such that a Hamilton vector field condition holds. To see what this is,
we consider a subset of continuous functions on $\dT X$, namely
the class of $\pi$-invariant $\Cinf$ functions on $T^* X$.
$\pi$-invariance means that if $\zeta,\zeta'\in T^*X$ and $\pi(\zeta)
=\pi(\zeta')$ then $f(\zeta)=f(\zeta')$. If $f$ is $\pi$-invariant then
it induces a function $f_\pi$ on $\dT X$ by $f_\pi(q)=f(\zeta)$ if
$q=\pi(\zeta)$. Moreover, if $f$ is smooth (or indeed just continuous)
then $f_\pi$ is continuous by the definition of
the topology on $\dT X$.

Now, if $\gammat$ is a curve in a manifold,
one way to put that it is an integral curve
of a vector field $V$ is that for all smooth functions $f$,
$\frac{d}{ds}(f\circ\gammat)|_{s=s_0}=(Vf)(\gammat(s_0))$.
If $f$ is a smooth $\pi$-invariant function on $T^*X$, then $f$
defines a $\Cinf$ function on $T^*X_a$ for all $a$, so $H_{g_a}f$ makes
sense. Here $H_{g_a}$ is the Hamilton vector field of the metric function
$g_a$ on $T^*X_a$,
so explicitly, $H_{g_a}=2\xi_a\cdot\pa_{x_a}$. Now we would like
to say that along a generalized broken bicharacteristic $\gamma$,
$\frac{d}{ds}(f_\pi\circ\gamma)|_{s=s_0}$ should be given by $H_{g_b}f$
for some $b$ and some $\zeta$ with $\pi(\zeta)=\gamma(s_0)$. The
problem is that there are many such points $\zeta$ and clusters $b$,
so this statement does not make any sense. However, we may replace the
derivative by the $\liminf$ of the difference quotients, i.e.\ by
$D_\pm h(s_0)=\liminf_{s\to s_0}\frac{h(s)-h(s_0)}{s-s_0}$, and demand
an inequality instead of the equality. That is, we may demand that
$D_\pm(f_\pi\circ\gamma)(s_0)$ may not be less than the worst possible
scenario as we run over all such $b$ and $\zeta$.
Thus, the condition for a
continuous map $\gamma:I\to\dChar(\lambda)$
to be a generalized broken bicharacteristic
is then that for any $s_0\in I$, if $\gamma(s_0)\in T^*X_{a,\reg}$
then
\begin{equation*}
D_\pm(f_\pi\circ\gamma)(s_0)\geq
\inf\{(H_{g_b}f)(\zeta):\ \zeta\in\Char_b(\lambda),\ \pi(\zeta)=\gamma(s_0),
\ X_a\subset X_b\}.
\end{equation*}
If the set of bound states is discrete, then such a curve $\gamma$ is
piecewise an integral curve of the Hamilton vector field of $g_b$
inside $\Char_b(\lambda)$, where $b$ may of course vary. In particular,
if there are no bound states in any proper subsystem, the picture
is very similar to wave propagation: the definition can be reduced
to the analogue of Lebeau's \cite{Lebeau:Propagation}.

The structure of the generalized broken
bicharacteristics, including the above claims, depends on having a large
supply of $\pi$-invariant functions. But these exist, since the pull-backs
of all functions on $X$ to $T^*X$ is $\pi$-invariant, so one can localize
in $X$ using smooth cutoffs. Moreover, near $X_{a,\reg}$, each component
of $\xi_a$ is $\pi$-invariant, as is $\xi^a\cdot x^a$. Note that
the generalized broken bicharacteristics depend on $V$, but only
via the characteristic set $\dChar(\lambda)$, i.e.\ only via
the bound states of the subsystem Hamiltonians.

There is also a wave front set associated to many-body scattering
which measures the microlocal decay of tempered distributions at infinity. For
a tempered distribution $u$, $\WFSc(u)$ is a closed conic subset of
$\dT X$. Apart from $u$, it depends on $\calX$, since $\dT X$ depends
on $\calX$, but we suppress this in the notation, and write
\begin{equation*}
\WFSc(u)=\WF_{\Scl,\calX}(u).
\end{equation*}
Its definition is slightly complicated, and I only refer to
\cite{Vasy:Propagation-Many} for the general definition,
which uses the structure of the pseudo-differential
algebra, in particular the operator-valued nature of symbols at
infinity. However, for generalized eigenfunctions of many-body Hamiltonians
it is simple. Namely, suppose that $(H-\lambda)u\in\Sch(X)$, where $\Sch(X)$
is the space of Schwartz functions. For $\bar x=\bar x_a\in X_{a,\reg}$ and
$\bar\xi_a\in X_a^*$ we say that $(\bar x_a,\bar\xi_a)\nin\WFSc(u)$
if there exists $\phi\in\Cinf_c(X_a^*)$ such that $\phi(\bar\xi_a)\neq 0$
and $\Frinv\phi\Fr u$ is rapidly decreasing in an open cone in $X$
around $\bar x_a$. Two examples are:
\begin{equation*}\begin{split}
\WFSc(e^{ix\cdot \xi_0})=\pi(\{(x,\xi_0):\ x\neq 0\}),\ \xi_0\in X_0^*,\\
\WFSc(e^{i\alpha|x|})=\pi(\{(x,\frac{\alpha x}{|x|}):\ x\neq 0\}),\ \alpha\in\Real.
\end{split}\end{equation*}
More generally, if $v$ is a symbol of any order on $X_0$, say
$v\in S^k(X_0)$, and $\phi\in\Cinf(X_0)$ is homogeneous degree $1$ for
$|x|>1$, then
\begin{equation*}
\WFSc(e^{i\phi(x)}v(x))\subset\pi(\Graph d\phi)
=\pi(\{(x,(d\phi)(x)): |x|> 1\}).
\end{equation*}
The condition $|x|>1$ is due to the requirement of the homogeneity
of $\phi$ only for $|x|>1$; technically we should add a subset
of $|x|\leq 1$ to the right hand side to make it conic.
The theorem on the propagation of singularities is the following.

\begin{thm}\label{thm:prop-sing}
Suppose that $\lambda\in\Real$ and $H$ is a many-body Hamiltonian.
If $u\in\Sch'(\Rn)$ and $(H-\lambda)u\in\Sch(\Rn)$ then
$\WFSc(u)\subset\dChar(\lambda)$ (microlocal elliptic regularity).
Moreover, $\WFSc(u)$ is a
union of maximally extended generalized broken bicharacteristics of
$H-\lambda$ (propagation of singularities).
\end{thm}

\begin{rem*}
In the time dependent version, one considers tempered distributional
solutions of $(D_t+H)u=0$ on $X_0\times\Real_t$. Then $D_t+H$ still
has the structure of a many-body Hamiltonian, with $D_t+\Delta$ in place of
$\Delta$, with collision planes given by $X_a\times\Real$, with $\{0\}$
added for the sake of completeness. Thus,
$t$ is always a variable along the collision planes, so in  particular,
its dual variable $\tau$, is $\pi$-invariant.
Moreover, $\Char_a(\lambda)$
is replaced by
\begin{equation*}
\Char_a=\{(x_a,t,\xi_a,\tau)\in T^*(X_a\times\Real):\ -\tau-g_a(\xi_a)\in
\pspec(H^a)\}\subset T^*X_a,
\end{equation*}
so effectively $-\tau$ plays the role of the energy $\lambda$.
Generalized broken bicharacteristics
can be defined as before with $H_{g_b}$ replaced
by $H_{\tau+g_b}=\pa_t+2\xi_b\cdot\pa_{x_b}$.
The main additional issue is that they can
only be expected to give a good description of propagation at finite
energies since $D_t+H$ is not elliptic in the usual sense.
So the analogue of Theorem~\ref{thm:prop-sing} is that
if $u\in\Sch'(\Rn\times\Real)$, $(D_t+H)u=0$, and $u=\psi(H)u$ for
some $\psi\in\Cinf_c(\Real)$, then $\WFSc(u)$ is a subset of the
characteristic set, and in fact $\WFSc(u)$ is a union of
maximally extended generalized broken bicharacteristics of $D_t+H$
inside it. The proof of this statement only requires simple modifications
of the proof of the theorem.
\end{rem*}

The interpretation of the theorem is much analogous to that for the
wave equation. However, there is a difference which also occurs
in the traditional microlocal setting for more general operators (i.e.\ for
operators other than the wave operator), see
\cite{Gullemin-Schaeffer:Fuchsian}. Namely, the orbits of the
$\Real^+$-action may be bicharacteristics, and then the statement of
the theorem is empty at the points lying on these orbits since the
wave front set is {\em a priori} conic. This happens for $(x_a,\xi_a)\in
T^*X_{a,\reg}\subset\dT X$ if and only if there exists some cluster $b$
with $X_b\supset X_a$, and $\zeta\in\Char_b(\lambda)$ such that $H_{g_b}$
at $\zeta$ is tangent to the orbits of the $\Real^+$-action. This happens,
in turn, if and only if $\xi_a$ is parallel to $x_a$ and $\lambda-|\xi_a|^2
\in\pspec(H^b)$. Such points are called radial points, and their
collection is denoted by
\begin{equation*}\begin{split}
\calR(\lambda)=\bigcup_{a\neq 1}\{(x_a,\xi_a)\in T^*X_{a,\reg}:
\ &\exists c\in\Real,\ \xi_a=cx_a,\\
&\exists b,\ X_b\supset X_a
,\ \lambda-|\xi_a|^2\in\pspec(H^b)\}.
\end{split}\end{equation*}
As we discuss in Section~\ref{sec:ac}, $\calR(\lambda)$ plays an important
role in asymptotic completeness. In many-body scattering it appeared
in the work of Sigal and Soffer \cite{Sigal-Soffer:N} and was called
`propagation set' because in the time-dependent picture this is where
particles end up as time goes to infinity. (In the stationary semiclassical
picture, this is where non-trapped
classical trajectories starting in a compact region end up.) It is
thus unfortunate, in terms of terminology,
that this is also the region where there is no real
principal type propagation.

\begin{rem*}
In the time-dependent problem, the set of radial points is
\begin{equation*}\begin{split}
\calR=\bigcup_a\{(x_a,t,\xi_a,\tau)\in T^*(X_{a,\reg}\times\Real):
&\ x_a=2t\xi_a,\\
&\exists b,\ X_b\supset X_a,
\ -\tau-|\xi_a|^2\in\pspec(H^b)\}.
\end{split}\end{equation*}
\end{rem*}

In terms of radial points, the difference between threshold energies
$\lambda\in\Lambda=\Lambda_0$ and non-threshold energies is that if
$\lambda\in\Lambda$, then there are {\em constant} generalized broken
bicharacteristics, i.e.\ bicharacteristics whose image is a single point.
Namely, if $x_a\in X_{a,\reg}$ and
$\lambda\in\pspec(H^a)$, then
\begin{equation*}
(x_a,0)\in \Char_a(\lambda)\cap T^*X_{a,\reg}\subset \dChar(\lambda),
\end{equation*}
and $H_{g_a}=2\xi_a\cdot\pa_{x_a}$ vanishes there, so $(x_a,0)$ is indeed
the image of a constant bicharacteristic. While this does not make
any difference for the propagation of singularities, it does for
the related limiting absorption principle, which in this generality
is due to Perry, Sigal and Simon \cite{Perry-Sigal-Simon:Spectral}.

\begin{thm}
If $\lambda\nin\Lambda$, then the limits
$R(\lambda\pm i0)=(H-(\lambda\pm i0))^{-1}$ exists as bounded operators
between $L^2_s$ and $H^2_{-s}$ for $s>\frac{1}{2}$. Here $H^m_l$ is
the weighted Sobolev space $\langle x\rangle^{-l}H^m(\Rn)$, $L^2_s=H^0_s$.
\end{thm}

In fact, the proofs of the limiting absorption principle and the propagation
of singularities are related. Indeed, the statement on propagation of
singularities can be strengthened for $R(\lambda+ i0)f$, $f\in\Sch(\Rn)$,
by saying that $\WFSc(R(\lambda+ i0)f)$ is not only a union of maximally
extended generalized broken bicharacteristics, as follows from
Theorem~\ref{thm:prop-sing}, but in fact it is a union of
generalized broken bicharacteristics $\gamma:\Real_s\to\dT X$
which go to
\begin{equation*}
\calR_+(\lambda)=\calR(\lambda)\cap
\bigcup_{a\neq 1}\{(x_a,\xi_a)\in T^*X_{a,\reg}:\ x_a\cdot\xi_a>0\}
\end{equation*}
as $s\to -\infty$. That is, the singularities at $\calR_+(\lambda)$ (where
the statement of Theorem~\ref{thm:prop-sing} is empty) can only leave
$\calR_+(\lambda)$ in the {\em forward} direction. The limiting
absorption principle is thus strengthened to:

\begin{thm}\label{thm:lim-abs-WF}
If $\lambda\nin\Lambda$, then for $f\in\Sch(\Rn)$,
$\WFSc(R(\lambda+ i0)f)$ is a subset of the image of $\calR_+(\lambda)$
under the {\em forward} generalized broken bicharacteristic relation.
A similar statement holds for $R(\lambda-i0)f$ with $\calR_+(\lambda)$ replaced
by
\begin{equation*}
\calR_-(\lambda)=\calR(\lambda)\cap
\bigcup_{a\neq 1}\{(x_a,\xi_a)\in T^*X_{a,\reg}:\ x_a\cdot\xi_a<0\},
\end{equation*}
and the forward relation by the backward relation.

In fact, if $u\in\Sch'(\Rn)$ and $\WFSc(u)$ is disjoint from the image of
$\calR_-(\lambda)$ under the backward generalized broken bicharacteristic relation,
then $R(\lambda+i0)u$ is defined by duality and $\WFSc(R(\lambda+i0)u)$
is a subset of the image of $\calR_+(\lambda)\cup\WFSc(u)$ under the
forward relation.
\end{thm}

\begin{rem*}
$\lambda\nin\Lambda$ can be also characterized by $\calR(\lambda)=\calR_+(\lambda)
\cup \calR_-(\lambda)$, i.e.\ that $x_a\cdot\xi_a$ never vanishes on $\calR(\lambda)
\cap T^*X_{a,\reg}$ for any $a$.

In the time-dependent setting, $x_a=2t\xi_a$ on $\calR$, so $x\cdot\xi_a
=2t$. So $\calR_+$, defined in $\calR$
by $x_a\cdot\xi_a>0$, is the subset of $\calR$
where $t>0$. Hence the `outgoing' terminology for $R(\lambda+i0)$ and
`incoming' for $R(\lambda-i0)$. In fact, the solution of $(D_t+H)u=0$
with $u|_{t=0}=\phi$, $\phi\in\Sch(X_0)$, say, is $u(.,t)=e^{-iHt}\phi$.
The time-dependent propagation of singularities shows that
$\WFSc(u)$ is a subset of the union of the image of $\calR_+$ under the forward
broken bicharacteristic relation and the image of $\calR_-$ under
the backward bicharacteristic relation.
Using the spectral measure and Stone's theorem,
\begin{equation*}
u(.,t)=\frac{1}{2\pi i}
\int_{\Real}e^{-i\lambda t}(R(\lambda+i0)-R(\lambda-i0))\phi\,d\lambda
\end{equation*}
Fixing some $\psi\in\Cinf_c(\Real)$, for $\phi$ in the range of $\psi(H)$
we thus deduce that in $t>0$, $\WFSc(u)$ arises from the
$R(\lambda+i0)$ term, and in $t<0$ from the $R(\lambda-i0)$ term.
So the time-dependent and stationary settings are very close: the only
difference is that in the latter, $\lambda$ is a parameter, while
in the former, it is a variable, $\lambda=-\tau$.
\end{rem*}

Again, one can make more precise propagation statements in some circumstances,
such as three-body scattering, where the precise nature of the singularities
can be analyzed, see \cite{Hassell:Plane, Vasy:Structure}.
Here we only state the stronger implication for the
structure of the scattering matrices, which we proceed to analyze.

\section{Scattering matrices}

Physically, the scattering matrices relate incoming and outgoing data in
an experiment. In the time independent framework (where $-\lambda$ is
the dual variable of time), for short-range potentials
an incoming wave of energy $\lambda$
in channel $\alpha$ (a channel is the choice of a cluster $a$ and
an $L^2$-eigenfunction $\psi_\alpha$ of $H^a$ of energy $\ep_\alpha$)
takes the following form in $|x|>1$:
\begin{equation*}
u_{\alpha,-}=e^{-i\sqrt{\lambda-\ep_\alpha}|x_a|}|x_a|^{-\frac{\dim X_a-1}{2}}
g_{\alpha,-}(\frac{x_a}{|x_a|})\psi_\alpha(x^a)+u'_-
\end{equation*}
Similarly, an outgoing wave has the form
\begin{equation*}
u_{\alpha,+}=e^{i\sqrt{\lambda-\ep_\alpha}|x_a|}|x_a|^{-\frac{\dim X_a-1}{2}}
g_{\alpha,+}(\frac{x_a}{|x_a|})\psi_\alpha(x^a)+u'_+
\end{equation*}
i.e.\ the sign of the phase has changed. Here $g_{\alpha,\pm}$ may
be taken e.g.\ $L^2$ functions on $S_a$, the unit sphere in $X_a$, or
ideally, at least one of them may be taken $\Cinf$. In either case,
$u'_\pm$ are `lower order terms', namely they
must be in $L^2_{-1/2}$. (Note that
$\langle x_a\rangle^{-\frac{\dim X_a-1}{2}}\in L^2_s(X_a)$
for $s<-1/2$ but not for $s=-1/2$.)
In fact, for $g_{\alpha,\pm}\in\Cinf_c(S_{a,\reg})$
we may take them to be of the form
$e^{-i\sqrt{\lambda-\ep_\alpha}|x_a|}|x_a|^{-\frac{\dim X_a+1}{2}}v$
where $v$ is a $0$th order symbol, with $S_{a,\reg}$ denoting
$X_{a,\reg}\cap S_0$.

One can now produce tempered distribution with given incoming, or alternatively
of given outgoing, asymptotics.
A typical example is of the form
\begin{equation}\label{eq:gen-ef}
P_{\alpha,+}(\lambda)g_{\alpha,-}
=u_{\alpha,-}-(H-(\lambda+i0))^{-1}((H-\lambda)u_{\alpha,-});
\end{equation}
here the lower order terms may be dropped from $u_{\alpha,-}$
without affecting $u=P_{\alpha,+}(\lambda)g_{\alpha,-}$
and $g_{\alpha,-}$ can be specified to be any smooth function on $S_a$.
In general, even if the incoming data are in a single channel $\alpha$,
as in \eqref{eq:gen-ef},
the corresponding generalized eigenfunction $u$ of $H$ will have outgoing waves
in all channels. The S-matrix $S_{\alpha\beta}(\lambda)$ picks out
the component in channel $\beta$ by projection in a certain sense,
see \cite{Vasy:Scattering}.
Thus, $S_{\alpha\beta}(\lambda)$ maps functions on $S_a$, the unit
sphere in $X_a$, to functions on $S_b$, by
\begin{equation*}
S_{\alpha\beta}(\lambda)g_{\alpha,-}=g_{\beta,+}
\end{equation*}
for $u$ as in \eqref{eq:gen-ef}.
For example, the free-to-free (i.e.\ $N$-cluster to $N$-cluster in
$N$-body scattering) S-matrix $S_{00}(\lambda)$ maps functions on $S_0$,
the unit sphere in $\Rn$, to functions on $S_0$, more precisely
$S_{00}(\lambda):L^2(S_0)\to L^2(S_0)$ is bounded.

More precisely, let $T_+$ be a pseudodifferential operator that is
identically $1$ on the outgoing radial set and identically $0$ on the incoming
radial set; see the paragraph of \eqref{eq:T_+-constr}
for a precise statement. Then
\begin{equation}\label{eq:S-pairing}
S_{\alpha\beta}(\lambda)=\frac{1}{2i\sqrt{\lambda-\ep_\beta}}
((H-\lambda)T_+P_{\beta,-}(\lambda))^*P_{\alpha,+}(\lambda),
\end{equation}
i.e.\ for any $g\in\Cinf(S_{a,\reg})$, $h\in\Cinf(S_{b,\reg})$,
\begin{equation*}
\langle h,S_{\alpha\beta}(\lambda)g\rangle
=\langle(H-\lambda)T_+P_{\beta,-}(\lambda)h,\frac{1}{2i\sqrt{\lambda-\ep_\beta}}
P_{\alpha,+}(\lambda)g\rangle.
\end{equation*}
This is equivalent to the usual wave operator definition in the time-dependent
setting, see \cite{Vasy:Scattering}. An immediate consequence of the
propagation of singularities and the definition of the scattering matrices
is the following:

\begin{thm}
The wave front relation of $S_{\alpha\beta}(\lambda)$ is given by the broken
bicharacteristic relation. In particular, if no proper subsystem of $H$
has bound states, the wave front relation of $S_{00}(\lambda)$ is given
by the broken geodesic flow on $S_0$ at distance $\pi$.
\end{thm}

While typically broken bicharacteristics can be continued in many ways
when they hit a collision plane, it is important to keep in mind that
under suitable assumptions (which rule out geometric complications)
the broken
bicharacteristic relation is Lagrangian, hence its dimension is the
same as if there were no collision planes. The reason is that only
a low dimensional family of broken
bicharacteristics hits any specified collision plane, with the dimension
of the possible continuations of each of these these bicharacteristics
compensating to yield the correct dimension for Lagrangian submanifolds.

The significant improvement in the three-body case, as shown by Hassell and
the author
\cite{Vasy:Structure, Hassell:Plane}, is that one can pinpoint not
only the location of the singularities, but also their precise form.
This theorem was motivated by the geometric result of Melrose and
Zworski \cite{RBMZw}, showing that the scattering matrix on asymptotically
Euclidean manifolds is a Fourier integral operator.

\begin{thm}
Suppose that $H$ is a three-body Hamiltonian and
the $V_a$ are Schwartz on $X^a$ for all $a$. Then
$S_{00}(\lambda)$ is a finite sum of Fourier integral operators (FIOs)
associated to the broken geodesic relation on $S_0$ to distance $\pi$.
Its canonical relation corresponds to the various collision patterns.
The principal symbol of the term corresponding to a single collision at
$X_a$ is given by, and in turn determines,
the 2-body S-matrix of $H^a$ at energies $\lambda'
\in(0,\lambda)$.
\end{thm}

\begin{rem*}
This result presumably extends to short range symbolic potentials, using the
same methods, though it is technically more complicated to write down the
argument in that case, and it has not been done. In fact, it should
also extend to the $N$-body problem, provided that there are no bound
states in any proper subsystem. Some assumption on the bound states
is necessary, for otherwise the generalize broken bicharacteristic relation
can become fairly complicated, see \cite{Vasy:Bound-States}.
The reason why one does not need any
assumption on bound states in three-body scattering is that
for any 2-cluster $a$, $\Char_a(\lambda)
\cap\pi_{0a}(\Char_0(\lambda))$ is either empty (if $0$ is not an
eigenvalue of $H^a$) or consists of the boundary of
$\pi_{0a}(\Char_0(\lambda))$. In the former case there is no interaction
(modulo smoothing terms) between the $0$-cluster and the $a$-cluster
dynamics, while in the latter case in the only place they interact, the
two dynamics give the same propagation.

It should also be noted that the normalization of $S_{\alpha\beta}(\lambda)$ is
not the standard one in many-body scattering
(which is based on wave operators), but rather follows
the geometric conventions \cite{RBMSpec}. The difference is that
in the wave operator approach, free motion is factored out, so the free
scattering matrix is the identity operator. On the other hand, in the
geometric approach we describe the asymptotics of generalized eigenfunctions,
or alternatively of the Schr\"odinger equation. Since free particles move
to infinity in the opposite direction from which they came, it is
reasonable that the two should differ by (a constant multiple of) pull-back by
the antipodal map, and this is indeed the case, see \cite{Vasy:Scattering}.
The distance $\pi$ propagation along (not broken!) geodesics on the sphere
indeed takes particles to the antipodal point.
\end{rem*}

An immediate corollary, when combined with two-body results
(e.g.\ analyticity of the S-matrix in $\lambda'$ and the Born approximation)
is the following inverse result.

\begin{cor*}
If the $V_a$ decay exponentially and $\dim X_a\geq 2$ for all $a$
then $S_{00}(\lambda)$ for a single value of $\lambda$ determines all
interactions.
\end{cor*}

This result is analogous to the recovery of cracks in a material
by directing sound waves at it and observing the singularities
of the reflected waves, except the
last step which uses two-body results to get the potentials from the
two-body S-matrices.

The other extremal scattering matrices are the 2-cluster to 2-cluster ones,
and they describe the physically most interesting events. Indeed, it is hard to
make more than two particles collide in an accelerator, so the initial state
in a physical experiment tends to be a 2-cluster. The following result
is due to Skibsted \cite{Skibsted:Smoothness}, and it also follows
from the propagation of
singularities and the definition of the S-matrices.

\begin{thm}
Let $\alpha$ and $\beta$ be two-clusters, and suppose that either
$\ep_\alpha\in\spec_d(H^a)$ and $\ep_\beta\in\spec_d(H^b)$, or
$V_c$ is Schwartz for all $c$.
Then the two-cluster to two-cluster S-matrix $S_{\alpha\beta}(\lambda)$
has $\Cinf$ Schwartz kernel,
except if $\alpha=\beta$ in which case the Scwartz kernel of
$S_{\alpha\alpha}(\lambda)$ is conormal to the graph of the antipodal map
on $S_a$,
corresponding to free motion.
\end{thm}

Thus, principal symbol calculations do not help in this inverse problem.
Note that if $H$ is a 3-body Hamiltonian, then
$\ep_\alpha\in\spec_d(H^a)$ and $\ep_\beta\in\spec_d(H^b)$ holds for
any non-threshold bound state energies.
The new result, in a joint project with Gunther Uhlmann,
is the following \cite{Uhlmann-Vasy:Low}.

\begin{thm}
Suppose that $H$ is a 3-body Hamiltonian, $a$ is a 2-cluster,
$\alpha$ is a channel of energy $\ep_\alpha<0$, $V_a$
is a symbol of negative order (i.e.\ may be long range). For any
$\mu>\dim X_a$ there
exists $\delta>0$ such that the following holds.

Suppose that $\sup|(1+|x^b|)^\mu V_b(x^b)|<\delta$ for all
$b\neq a$. Suppose also that $I\subset(\ep_\alpha,0)$ is a non-empty
open set, and let
\begin{equation*}
R=2\sqrt{\sup I-\ep_\alpha}.
\end{equation*}
Then $S_{\alpha'\alpha''}(\lambda)$ given for all $\lambda\in I$ and
for all bound states $\alpha',\alpha''$ of $H^a$ with $\ep_{\alpha'},
\ep_{\alpha''}<\sup I$, determines the
Fourier transform of the effective interaction $V_{\alpha,\eff}$
in the ball of radius $R$ centered at $0$.
\end{thm}

The effective interaction is the interaction that arises if we consider
the 3-body problem as a 2-body problem, i.e.\ if we regard the two
particles forming the cluster $a$ as a single particle. Mathematically,
this amounts to projecting to the state $\psi_\alpha$ in $X^a$ and
obtaing a new Hamiltonian $\Delta_{X_a}+V_{\alpha,\eff}$ on $X_a$.
Thus, the effective interaction is physically relevant. Moreover,
there is no hope for recovering anything better than $V_{\alpha,\eff}$
as shown by the high-energy inverse results of Enss and Weder
\cite{Enss-Weder:Geometrical, Enss-Weder:Inverse},
Novikov \cite{Novikov:N-body} and Wang \cite{Wang:High, Wang:Inverse}.

This theorem says that if the unknown interactions are small then the
effective interaction can be determined from the knowledge of
all S-matrices with incoming and outgoing data in the cluster $a$ in
the relevant energy range.
In fact, near-forward information suffices as in two-body scattering,
where this was observed recently by Novikov \cite{Novikov:Determination}.
Also, if one is willing to take small $R$ and
$\alpha$ is the ground state of $H^a$,
it suffices to know $S_{\alpha\alpha}(\lambda)$
to recover $\hat V_{\alpha,\eff}$ in a small ball.

In case $V_b$ decay exponentially on $X^b$ for all $b\neq a$, then
$V_{\alpha,\eff}$ decays exponentially on $X_a$, hence its Fourier
transform is analytic, so $V_{\alpha,\eff}$ itself can be recovered
from these S-matrices.

\begin{rem*}
It is clear from the proof in \cite{Uhlmann-Vasy:Low} that there is
a natural extension of this theorem to many-body scattering at low energies.
\end{rem*}

This result should extend to higher energies, i.e.\ $\sup I\leq 0$
is {\em not} expected to be essential. But it is hard to make
$R$ greater than $2\sqrt{-\ep_\alpha}$ even then.
The reason is that our method relies on the construction of exponential
solutions following Faddeev \cite{Faddeev:Inverse},
Calder\'on \cite{Calderon:Inverse}, Sylvester and Uhlmann
\cite{Sylvester-Uhlmann:Global} and Novikov and
Khenkin~\cite{Novikov-Khenkin:D-bar}, but in the three-body setting.
One thus allows complex momenta $\rho\in\Cx(X_a)$, the complexification
of $X_a$, and one wants to construct solutions of $(H-\lambda)u=0$ of
the form
\begin{equation*}
e^{i\rho\cdot x_a}(\psi_\alpha(x^a)+v),
\end{equation*}
where $v=v_\rho$ is supposed to be `small' in the sense that it goes
to $0$ as $\rho\to\infty$ in an appropriate fashion. Note that with $v=0$
these complex plane waves solve $(H_a-\lambda)u=0$ with
\begin{equation}\label{eq:cx-energy}
\lambda=\rho\cdot\rho+\ep_\alpha;
\end{equation}
this expresses that the total energy $\lambda$ is the sum of the kinetic
energy, $\rho\cdot\rho$, and the potential energy $\ep_\alpha$.

To construct $u$, we need to find $v$, and its study reduces to that of the
conjugated Hamiltonian
\begin{equation*}
e^{-i\rho\cdot x_a}(H-\lambda)e^{i\rho\cdot x_a}
=H^a+\Delta_{X_a}+2\rho\cdot D_{X_a}+I_a-\ep_\alpha
\end{equation*}
with $\rho\in\Cx(X_a)$ the complex frequency.
Here we used \eqref{eq:cx-energy}. Now, $I_a$ is considered as a perturbation
(this is the reason for the smallness assumption in the theorem),
so we really study the model operator,
\begin{equation*}
H^a+\Delta_{X_a}+2\rho\cdot D_{X_a}-\ep_\alpha.
\end{equation*}
Taking the Fourier transform in the $X_a$ variables, one obtains
\begin{equation*}
H^a+|\xi_a|^2+2\rho\cdot \xi_a-\ep_\alpha.
\end{equation*}
Writing $\rho=z\nu+\rho_\perp$ with $|\nu|=1$,
$\rho_\perp\cdot\nu=0$, $\rho,\nu$ real, $z\in\Cx$, this operator becomes
\begin{equation*}
H^a+|\xi_a|^2+2\rho_\perp\cdot \xi_a+2z\nu\cdot \xi_a-\ep_\alpha
=H^a+(\xi_a+\rho_\perp)^2+2z\nu\cdot\xi_a-|\rho_\perp|^2-\ep_\alpha.
\end{equation*}
If $\rho$ is not real, then neither is $z$, so this operator is invertible
if $\nu\cdot\xi_a\neq 0$ since $H^a$ is self-adjoint. On the other
hand, if $\nu\cdot\xi_a=0$, this operator becomes
\begin{equation*}
H^a+(\xi_a+\rho_\perp)^2-|\rho_\perp|^2-\ep_\alpha,
\end{equation*}
i.e.\ its
invertibility properties correspond to the behavior of the boundary
values of the resolvent of $H^a$ at the real axis. If
$|\rho_\perp|^2+\ep_\alpha<0$, i.e.\ if $|\rho_\perp|<\sqrt{-\ep_\alpha}$,
then the spectral parameter $|\rho_\perp|^2+\ep_\alpha-(\xi_a+\rho_\perp)^2$
is negative, so only the bound states of $H^a$ contribute to the
characteristic variety, i.e.\ the two-cluster $a$ may not break up.
On the other hand, if $|\rho|\geq \sqrt{-\ep_\alpha}$, such a break-up
is possible {\em even if} $\lambda<0$, i.e.\ where the break up may not
happen for {\em real} frequencies.
The break-up greatly influences analyticity properties, hence one cannot
easily use large $\rho_\perp$.
On the other hand, one needs such large
$\rho_\perp$ to recover $V_{\alpha,\eff}$ on larger balls, hence the
limitation in the theorem. This also suggests that the fixed energy problem
would be hard, since then one always needs to let $\rho_\perp\to\infty$
to keep $\rho\cdot\rho=|\rho_\perp|^2+z^2$ fixed and yet have $\rho\to\infty$.

\section{Many-body scattering pseudo-differential operators}\label{sec:psdo}
I will present the calculus from the compactified point of view.
Both the one-step
polyhomogeneous (i.e.\ `classical') and the non-polyhomogeneous
calculus can be described in non-compact terms, i.e.\ directly on $X_0$,
but this is more complicated and less natural. Indeed, one of the
beauties of compactification is that it exactly captures the structure
of many-body Hamiltonians. {\em We warn the reader here that from
now on the Euclidean variable is written as $z$, rather than $x$ in
the preceeding sections, for compatibility with previous papers
espousing this approach, such as \cite{RBMSpec, RBMGeo}.}

To see how the compactification should go, recall first that a
classical symbol of order $0$ on $\Rn_z$ has an asymptotic expansion
\begin{equation*}
a(r\omega)\sim\sum_{j=0}^\infty r^{-j}a_j(\omega),\ a_j\in\Cinf(\Sn),
\end{equation*}
in the polar coordinates $(r,\omega)$: $z=r\omega$. The meaning of such
an expansion is that, for any $k$,
the difference of $a$ and the sum of the first $k$
terms on the right hand side is a symbol of order $-k$.
This expansion is just a Taylor series at $r=\infty$, or rather at
`$r^{-1}=0$'. So we compactify $\Rn$ into a ball $\overline{\mathbb{B}^n}$
by adding points $(0,\omega)$, $\omega\in\Sn$, and making
$(r^{-1},\omega)=(x,\omega)$ coordinates near these points. The
resulting space is
called the radial compactification $\overline{\Rn}$ of $\Rn$.
Thus, a classical symbol of order $0$ is simply a smooth function
of $\overline{\Rn}$; the asymptotic expansion at infinity is its
Taylor series around the boundary, $x=0$.

This compactification, whose utility in this context was emphasized
by Melrose \cite{RBMSpec},
can also be realized as the closed unit upper hemisphere
via a modified stereographic projection. So let
\begin{equation*}
\SP:\Rn\to\Snp,\ \SP(z)=(\frac{1}{\langle z\rangle},\frac{z}{\langle
z \rangle}),\ \langle z \rangle=(1+|z|^2)^{1/2},\ z\in\Rn.
\end{equation*}
Then $n$ of the $n+1$ variables $(\frac{1}{\langle z\rangle},\frac{z}{\langle
z \rangle})$ give local coordinates on various regions of $\Snp$.
In particular, in coordinate patches near the equator, which is $\pa\Snp$,
$\frac{1}{\langle z\rangle}$ (or indeed $x=|z|^{-1}$) and
$n-1$ of $\frac{z_j}{\langle z\rangle}$ (or indeed $\omega_j=\frac{z_j}{|z|}$)
can be taken as coordinates, showing that $\Snp$ can be identified
with the radial compactification $\overline{\Rn}$. A slightly
modified version of $x$ (it needs to be smoothed at $z=0$, where `$x=\infty$'),
or $\langle z\rangle^{-1}$, can be taken as a boundary defining function.
We will usually write $x$ for this, so $x=|z|^{-1}$ for $|z|\geq 1$, say.
(A boundary defining function is a non-negative function whose zero set
is exactly the boundary, and whose differential does not vanish there.)

\begin{figure}
\begin{center}
\mbox{\epsfig{file=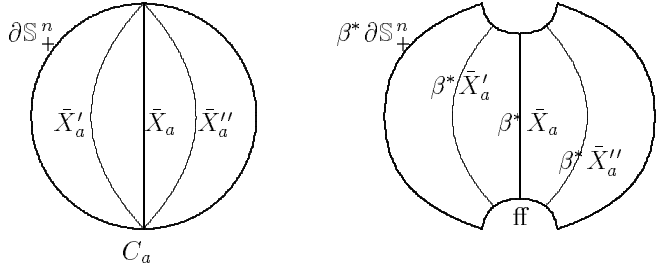}}
\end{center}
\caption{Translates of $X_a$ on $[\Xb_0;C_a]$.}
\label{fig:ressp}
\end{figure}

How can we adapt this to many-body scattering? Let $\Xb_a$ denote the closure
of $X_a$ in the compactification $\overline{\Rn}$ of $\Rn$,
and let $C_a=\pa\Xb_a\subset\pa\Snp=C_0$. The closure
of any translate of $X_a$ intersects $C_0$ in the same submanifold (a sphere)
as $X_a$ itself. Indeed, writing the coordinates as $(z_a,z^a)$ on
$X_0=X_a\oplus X^a$, local coordinates near $C_a$ are given by
$Z^a=\frac{z^a}{|z|}$, $|z|^{-1}$ and $\dim X_a-1$ of $\frac{(z_a)_j}{|z|}$.
Thus, $Z^a\to 0$ as $x\to 0$ along any translate, since $z^a$ is constant
along these. So $V_a$ is not even
continuous on $\Xb_0$, as it takes different values on the different translates
of $X_a$. However, it {\em is} a negative order symbol (in particular
continuous with boundary value $0$) on $\Xb_0\setminus C_a$, if $V_a$
is such on $X^a$; see Figure~\ref{fig:ressp}.

\begin{figure}
\begin{center}
\mbox{\epsfig{file=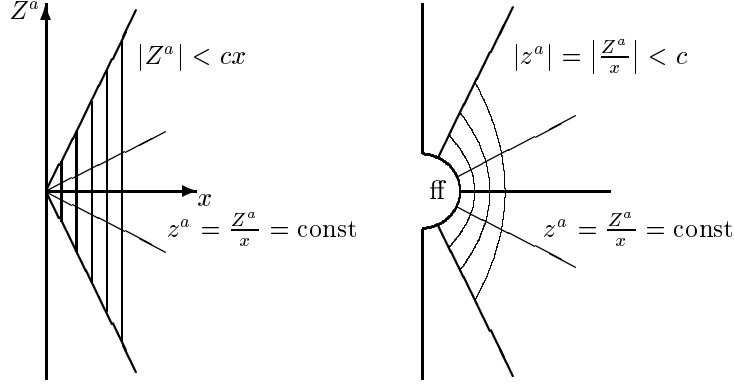}}
\end{center}
\caption{The blow up of $C_a$, given by $Z^a=0$, $x=0$.}
\label{fig:blowup}
\end{figure}

So the compactification works for $V_a$, except at $C_a$. To remedy this,
we blow up $C_a$. This is an invariant way of introducing
polar coordinates about it
(i.e.\ projective coordinates in various charts). That is, curves
approaching $C_a$ from various normal directions will correspond to different
points on the blown-up space $[\Xb_0;C_a]$. Since $C_a$ is
given by $x=0$, $Z^a=0$, in local coordinates, this means concretely
that the components of $Z^a/x$ become coordinate functions on the
part of $[\Xb_0;C_a]$ where this quotient is finite. (For the sake of
completeness, a complete set of coordinates in this region is given by
$x$, the components of $Z^a$ as well as the $\dim X_a-1$ coordinates on
the sphere $y_a=\frac{z_a}{|z_a|}$; see Figure~\ref{fig:blowup}.)
But $Z^a/x=z^a$, so
it is now easy to see that for classical symbols $V_a$ on $X^a$
(of negative integer
order), $V_a$ is a $\Cinf$ function on $[\Xb_0;C_a]$.

In general, there are many collision planes, and we blow them up
recursively, starting with ones of the largest codimension, to get
$[\Xb_0;\calC]$,
\begin{equation*}
\calC=\{C_a:\ X_a\in\calX,\ a\neq 1\}.
\end{equation*}
We refer
to \cite{Vasy:Propagation-Many} for details.

There is no reason at all to take $\Xb_0$ as the space we start with.
Given any compact
manifold with boundary, $\Xb$, and a cleanly intersecting family
of closed embedded submanifolds $\calC$ of $\pa\Xb$, we can define
$[\Xb;\calC]$ analogously. For instance, one can start with
$\Xb=\overline{\Real^n}\times\sphere^k$.
The space $[\Xb;\calC]$ is equipped with boundary fibrations given
by the blow-down maps, see \cite{Mazzeo-Melrose:Fibred} for a simpler
case where these first appeared explicitly.

Having described the configuration space, we turn to differential
operators. $X_0$ has a nice algebra of differential operators,
consisting of operators with symbolic coefficients:
$\sum_{|\alpha|\leq m}a_\alpha(z)D_z^\alpha$, $a_\alpha\in S^0(X_0)$.
We may require instead that $a_\alpha$ is `classical', i.e.\ that
$a_\alpha\in\Cinf(\Xb_0)$. The resulting algebras were denoted
$\Diffscc(\Xb_0)$ and $\Diffsc(\Xb_0)$ by Melrose; he called them
`scattering differential operators'.

This setup generalizes to the geometric set-up as follows. Let $(x,y)$,
$y=(y_1,\ldots,y_{n-1})$, be
local coordinates near $\pa\Xb_0$. Then the vector fields in
$\Diffsc(\Xb_0)$ are linear combinations of $x^2D_x$ and the $xD_{y_j}$
with coefficients in $\Cinf(\Xb_0)$, as can be seen easily by
an explicit calculation.

Now, if $\Xb$ is a manifold with boundary, $\Vb(\Xb)$ is the Lie
algebra of vector fields tangent to $\pa\Xb$, and $\Vsc(\Xb)=x\Vb(\Xb)$,
where $x$ is a defining function of $\pa\Xb$. $\Vsc(\Xb)$ is independent
of the choice of $x$. Then $\Vb(\Xb)$ is spanned by $x\pa_x$ and
$\pa_y$ over $\Cinf(\Xb)$, so $\Vsc(\Xb)$ is spanned by
$x^2\pa_x$ and $x\pa_y$ over $\Cinf(\Xb)$. By definition, these
generate $\Diffsc(\Xb)$. Also, $\Vsc(\Xb)$ is the set of all smooth
sections of a vector bundle over $\Xb$, this is denoted by
$\Tsc \Xb$. Its dual bundle is the scattering cotangent bundle,
denoted by $\sct \Xb$. In the Euclidean setting,
\begin{equation*}
\Tsc\Xb_0=\Xb_0\times X_0,\ \sct\Xb_0=\Xb_0\times X_0^*.
\end{equation*}

The way to generalize this differential operator algebra to one that
includes many-body potentials is to allow singular coefficients
$a_\alpha\in\Cinf([\Xb;\calC])$. Thus,
\begin{equation*}
\DiffSc(\Xb;\calC)=\Cinf([\Xb;\calC])\otimes_{\Cinf(\Xb)}\Diffsc(\Xb).
\end{equation*}
In particular, if $H$ is a many-body Hamiltonian, with either potential
or higher order interactions, then $H\in\DiffSc^2(\Xb_0;\calC)$.

We let $\TSc[\Xb;\calC]=\beta^*\Tsc\Xb$ and $\Tsc^*[\Xb;\calC]
=\beta^*\Tsc^*\Xb$, where $\beta:[\Xb;\calC]
\to\Xb$ is the blow-down map, and we are pulling back the vector bundles
by it. Hence in the Euclidean setting,
\begin{equation*}
\TSc[\Xb_0;\calC]=[\Xb_0;\calC]\times X_0,
\ \TSc^*[\Xb_0;\calC]=[\Xb_0;\calC]\times X_0^*.
\end{equation*}
Now it is natural to define pseudo-differential operators using these
bundles. Although I restrict the discussion to the Euclidean setting,
the construction generalizes to any $\Xb$ via localization.

So we consider symbols
\begin{equation}\label{eq:symbol-class}
a\in\langle z\rangle^{-l}\langle \zeta\rangle^m\Cinf([\Xb_0;\calC]
\times \Xb_0^*),
\end{equation}
$X_0=\Rn$, where $\zeta$ is the dual variable of $z$,
i.e.\ the variable on $X_0^*$. Note that this means that $a$ is a
classical symbol of order $m$ in $\zeta$. As usual,
we define the Schwartz kernel of
the left quantization of $a$ by
\begin{equation}\label{eq:left-q}
A=q_L(a)=(2\pi)^{-n}\int_{\Rn} e^{i(z-z')\cdot\zeta}a(z,\zeta)\,d\zeta,
\end{equation}
understood as an oscillatory integral. In particular, for any $f\in
\Sch(\Rn)$,
\begin{equation*}
Af(z)=(2\pi)^{-n}\int_{\Rn}\int_{\Rn}
e^{i(z-z')\cdot\zeta}a(z,\zeta)f(z')\,d\zeta\,dz',
\end{equation*}
again understood as an oscillatory integral. We write
$A\in\PsiSc(\Xb_0;\calC)$ for this class of operators.

Note that $\langle z\rangle^l a\in S^m_\infty(X_0;X_0^*)$, H\"ormander's
uniform symbol space \cite[Section~18.1]{Hor}, so $A=\langle z\rangle^l
\tilde A$, $\tilde A\in\Psi^m_\infty(X_0)$, the uniform ps.d.o.-algebra
arising by quantizing $S^m_\infty(X_0;X_0^*)$ as in \eqref{eq:left-q}.
In particular, since the
mapping properties of $\Psi^m_\infty(X_0)$ between weighted Sobolev
spaces $H^{r,s}$ are well known, the corresponding properties of $A$
follow. Namely, $A:H^{r,s}\to H^{r-m,s+l}$ for all $r,s$, where
\begin{equation*}
H^{r,s}=\langle z\rangle^{-s}H^r=\langle f\in\Sch'(\Rn):\ \langle z\rangle
^s f\in H^r\}.
\end{equation*}

Now, $\PsiSc(\Xb_0;\calC)$ is a $*$-algebra, in particular is
closed under composition. Indeed,
since $\PsiSc(\Xb_0;\calC)\subset\Psi_\infty(X_0)$, and the
latter is closed under composition, it suffices to follow the usual
proof and make sure that the product is in $\PsiSc(\Xb_0;\calC)$,
rather than merely in $\Psi_\infty(X_0)$. Thus, the key fact is that
for any
\begin{equation*}
b\in\langle z\rangle^{-l}\langle \zeta\rangle^m\Cinf([\Xb_0;\calC]_z
\times[\Xb_0;\calC]_{z'}\times(\Xb_0^*)_\zeta)
\end{equation*}
there exists $a$ as in \eqref{eq:symbol-class} such that
the induced operators
\begin{equation}\label{eq:sym-quant}
B=(2\pi)^{-n}\int_{\Rn} e^{i(z-z')\cdot\zeta}b(z,z',\zeta)\,d\zeta,
\end{equation}
and $A$ as in \eqref{eq:left-q} are the same. The proof of this
claim is standard. Indeed, we can expand $b$ in Taylor series
in $z'$ around $z=z'$ to finite order $k$. The finite order terms depend
on $z'$ only via $(z'-z)^\alpha$, $|\alpha|\leq k$. We rewrite
$(z'-z)^\alpha e^{i(z-z')\cdot\zeta}$ as $(-1)^{|\alpha|}D_\zeta^\alpha
e^{i(z-z')\cdot\zeta}$, and integrate by parts with respect to $\zeta$.
Thus, the $\alpha$-term is the left quantization of
\begin{equation}\label{eq:prod-terms}
\frac{1}{\alpha!}D^\alpha_{z'}D^\alpha_\zeta b(z,z',\zeta)|_{z'=z},
\end{equation}
which is of the desired form, i.e.\ is in
$\langle z\rangle^{-l}\langle \zeta\rangle^m\Cinf([\Xb_0;\calC]
\times \Xb_0^*)$. In fact, the weight $\langle \zeta\rangle^m$
can be replaced by $\langle \zeta\rangle^{m-|\alpha|}$ due to the symbolic
properties of $b$ in $(X_0^*)_\zeta$, but no corresponding change
may be made for the $z$ weight.
Similarly, the remainder term is
of the form 
\begin{equation}\begin{split}\label{eq:comp-rem}
&K_k(z,z')=(2\pi)^{-n}\int_{\Rn} e^{i(z-z')\cdot\zeta}b_k(z,z',\zeta)\,d\zeta,\\
&b_k\in\langle z\rangle^{-l}\langle \zeta\rangle^{m-k-1}\Cinf([\Xb_0;\calC]_z
\times[\Xb_0;\calC]_{z'}\times(\Xb_0^*)_\zeta).
\end{split}\end{equation}
Now we can asymptotically sum the $b_\alpha$ to get a new symbol
\begin{equation*}
c\in\langle z\rangle^{-l}\langle \zeta\rangle^m
\Cinf([\Xb_0;\calC]\times\Xb_0^*).
\end{equation*}
Let $C$ be the left quantization of
$c$. Then $B-C$ is of the form \eqref{eq:comp-rem} for all $k$, with
$b_k$ replaced
by some $b_k'$ with the same properties. It is then straightforward
to show that the Schwartz kernel $K'$ of $B-C$ is $\Cinf$, decays
rapidly with all derivatives as $\langle z-z'\rangle\to\infty$, and
more precisely it is of the form
\begin{equation*}
K'\in\Cinf([\Xb_0;\calC]_z\times(\Xb_0)_{z-z'})
\end{equation*}
with infinite order vanishing at the boundary of the second factor.
Taking its Fourier transform $b'$ in $z-z'$, $K'$ is thus the left
quantization of $a=c+b'$, proving the claim, hence in turn that
$\PsiSc(\Xb_0;\calC)$ is closed under composition.

In the two-body setting, where $\calC=\emptyset$, there is a principal symbol
at infinity. Namely, if $A\in\PsiSc^{m,l}(\Xb)$, $A=q_L(a)$, then
$\sigma_{m,l}(A)$ is given by
the restriction of
$\langle z\rangle^l\langle \zeta\rangle^{-m}
a\in\Cinf(\Xb_0\times\Xb_0^*)$ to $\pa(\Xb_0\times\Xb_0^*)
=(\pa\Xb_0\times\Xb_0^*)\cup(\Xb_0\times\pa\Xb_0^*)$.
Of the two boundary hypersurfaces,
the restriction to $\Xb_0\times\pa\Xb_0^*$ yields the usual principal
symbol, while the restriction to $\pa\Xb_0\times\Xb_0^*$ is the principal
symbol at infinity. More precisely, if $l=0$, we can indeed define
the part of $\sigma_{m,0}(A)$ at infinity to be the restriction of $a$ to
$(\pa\Xb_0)\times\Xb_0^*$.
The principal symbol is multiplicative,
i.e.\ $\sigma_{m+m',l+l'}(AB)=\sigma_{m,l}(A)\sigma_{m',l'}(B)$.
Thus, $[A,B]\in\PsiSc^{m+m'-1,l+l'+1}(\Xb)$, and its principal symbol
is given by the Poisson bracket of their symbols, see
Section~\ref{sec:positive}.

Since in the many-body setting
we do not gain decay in $z$ in \eqref{eq:prod-terms},
we cannot expect to have
a commutative principal symbol at infinity, i.e.\ at $\pa[\Xb_0;\calC]$.
For $C\in\PsiSc^{m,0}(\Xb,\calC)$,
$y_a\in C_{a,\reg}$, $\zeta_a\in X_a^*$, we let
\begin{equation*}
\Ch_a(y_a,\zeta_a)=(2\pi)^{-\dim X^a}\int e^{i(z^a-(z')^a)\cdot\zeta^a}
c(y_a,z^a,\zeta)\,d\zeta\in\Sch'(X^a\times X^a)
\end{equation*}
be the operator valued principal symbol of $C$ at $(y_a,\zeta_a)$.
Thus, $\Ch_a(y_a,\zeta_a)$ is a tempered distribution on $X^a\times X^a$
(denoted by the variables $(z^a,(z^a)')$), and it is in fact a
many-body ps.d.o. itself: $\Ch_a(y_a,\zeta_)\in \PsiSc^{m,0}(\Xb^a,\calC^a)$
corresponding to the collision planes $X^a\cap X_b$, with $b$ satisfying
$X_b\supset X_a$.
We also call it the indicial operator of $C$ to make it clear we are
not talking about the standard principal symbol. We also write
$\Ch_a(z_a,\zeta_a)$ in the same setting, where we extend $\Ch_a(y_a,\zeta_a)$
to be homogeneous degree $0$ in $z_a$.
It can be easily seen to satisfy
\begin{equation*}
\Ah_a\Bh_a=\widehat{(AB)}_a,
\end{equation*}
where on the left hand side we compose the operators
$\Ah_a(z_a,\zeta_a)$ and $\Bh_a(z_a,\zeta_a)$.
Thus, multiplication of operators is only partially commutative, even to
top order. This can be observed already from $[D_{z_a},V_a]=0$,
hence certainly lower order at infinity, while $[D_{z^a},V_a]
\in\Cinf([\Xb_0;\calC])$ without {\em any} decay at $C_a$.

This observation has important implications for the positive commutator
estimates that we take up in the next section.
Namely, $H$ must commute to leading order
with the operators we want to microlocalize with.
This means that these operators $A$ must have $\Ah_a$ commute with
$\Hh_a$, and the most reasonable way of achieving this is to have
$\Ah_a$ be a scalar multiple of $\psi(\Hh_a)$, where e.g.\ $\psi\in
\Cinf_c(\Real)$. This multiple
defines a function on $\dT\Xb_0$; we want this to arise from a smooth
$\pi$-invariant function for our estimates. On the other hand,
$\psi(\Hh_a)$ provides localization at the characteristic set.

Here, however, I would like to talk about pseudo-differential constructions
first. Namely, if $\lambda\nin\Real$, or indeed if $\lambda\in
\Cx\setminus[\inf\Lambda,+\infty)$ then there exists a parametrix
$G(\lambda)\in\PsiSc^{-2,0}(\Xb_0;\calC)$ for $H-\lambda$, i.e.\ such
that
\begin{equation*}
(H-\lambda)G(\lambda)-\Id,
\ G(\lambda)(H-\lambda)-\Id\in\PsiSc^{-\infty,\infty}
(\Xb_0;\calC).
\end{equation*}
Then the parametrix identities show that
\begin{equation*}
\lambda\in\Cx\setminus\spec(H)\Rightarrow
(H-\lambda)^{-1}\in\PsiSc^{-2,0}(\Xb_0;\calC).
\end{equation*}
The parametrix
construction proceeds inductively
by constructing $(\Hh_a-\lambda)^{-1}$ in $\PsiSc^{-2,0}(\Xb_a,\calC^a)$
for every $a\neq 1$
and then combining these: there exists a $G_0(\lambda)\in\PsiSc^{-2,0}
(\Xb_0;\calC)$ with specified indicial operators $(\Hh_a-\lambda)^{-1}$,
hence satisfying
\begin{equation*}
(H-\lambda)G_0(\lambda)-\Id,
\ G_0(\lambda)(H-\lambda)-\Id\in\PsiSc^{-1,1}
(\Xb_0;\calC).
\end{equation*}
Then the standard Neumann series argument yields $G(\lambda)$.

The Helffer-Sj\"ostrand argument \cite{Helffer-Sjostrand:Schrodinger}
then shows that for any $\phi\in\Cinf_c(\Real)$,
\begin{equation*}
\phi(H)=\frac{-1}{2\pi i}\int_{\Cx}\overline\pa_\lambda\tilde\phi(\lambda)
(H-\lambda)^{-1}\,d\lambda\wedge d\bar\lambda,
\end{equation*}
where $\tilde\phi$ is an almost analytic extension of $\phi$:
$\tilde\phi\in\Cinf_c(\Cx)$, $|\overline\pa_\lambda\tilde\phi|\leq
C_k|\im\lambda|^k$ for all $k$. We can control $(H-\lambda)^{-1}$
in $\PsiSc^{-2,0}(\Xb,\calC)$ as $\lambda\to\Real$ with semi-norm estimates
$\calO(|\im\lambda|^{-j})$ ($j$ depends on the norm), so we conclude that
$\phi(H)\in\PsiSc^{-\infty,0}(\Xb;\calC)$.

We can now explain the precise specifications on $T_+$ in \eqref{eq:S-pairing}.
Namely, we require that on a neighborhood of $\calR_+(\lambda)$ in
$\dT X$, the indicial operators $\widehat{T_+}$ should equal
$\widehat{\phi(H)}$ for some $\phi\in\Cinf_c(\Real)$ identically $1$
near $\lambda$, and on a neighborhood of $\calR_-(\lambda)$ they should vanish.
Explicitly this can be arranged by taking any $\phi$ as above, and any
$\chi\in\Cinf(\Real)$ identically $1$ on $(\frac{\sqrt{\lambda}}{2},+\infty)$,
identically $0$ on $(-\infty,-\frac{\sqrt{\lambda}}{2})$. Then let
$T_+=\phi(H) q_R(\chi(\frac{z\cdot\zeta}{\langle z\rangle}))$, with
$q_R$ denoting the `right quantization' (i.e.\ where we take
$b=\chi(\frac{z'\cdot\zeta}{\langle z'\rangle})$ in \eqref{eq:sym-quant}).
Although $q_R(\chi(\frac{z\cdot\zeta}{\langle z\rangle}))$ is not
in $\PsiSc(\Xb;\calC)$,
due to the non-symbolic behavior of $b$ as $\zeta\to\infty$,
$T_+$ is, namely $T_+\in\PsiSc^{-\infty,0}(\Xb;\calC)$,
since $\phi(H)$ is smoothing: see \cite{Vasy:Propagation-Many}.
Moreover,
\begin{equation}\label{eq:T_+-constr}
\widehat{T_+}_a(z_a,\zeta_a)=\chi(\frac{z_a\cdot\zeta_a}{|z_a|})
\phi(\Hh_a(z_a,\zeta_a)),
\end{equation}
hence has the desired properties.

Our construction of $\phi(H)$
in fact shows that if all potentials are in $S^{-\rho}(X_a)$,
$\rho>0$, and $\chi_a\in\Cinf(\Xb_0)$ is supported away from $C_b$
such that $C_b\supset C_a$ does {\em not} hold, then
$\chi_a(\phi(H)-\phi(H_a))\in\PsiSc^{-\infty,\rho}(\Xb;\calC)$, hence
trace class if $\rho>\dim X_0$. In the three-body setting
this shows that
\begin{equation*}
\phi(H)-\phi(H_0)-\sum_{\#b=2}(\phi(H_b)-\phi(H_0))
\end{equation*}
is trace class. Indeed, near $C_a$ this can be written as
\begin{equation*}
(\phi(H)-\phi(H_a))-\sum_{\#b=2,\ b\neq a}(\phi(H_b)-\phi(H_0)),
\end{equation*}
and now all terms in parantheses are in $\PsiSc^{-\infty,\rho}(\Xb;\calC)$
near $C_a$. So we conclude, with a proof that shows much more,
a result of Buslaev and Merkureev:
\begin{equation*}
\sigma(\phi)=\tr((\phi(H)-\phi(H_0))-\sum_{\#b=2}(\phi(H_b)-\phi(H_0)))
\end{equation*}
defines a distribution $\sigma\in\dist(\Real)$. Writing $\sigma=\xi'$
defines the spectral shift function, up to a constant, which in
turn, in two-body scattering,
is the well-known generalization of the eigenvalue counting function
on compact manifolds.
These statements, as well as the following theorem, which is
joint work with Xue-Ping Wang \cite{Vasy-Wang:Smoothness}, generalize to
arbitrary many-body Hamiltonians (with short-range interactions as indicated).

\begin{thm}
Suppose $H$ is a three-body Hamiltonian with Schwartz interactions:
$V_a\in\Sch(X^a)$, and that all interactions are pair interactions
(i.e.\ $V_a\neq 0$ implies that $a$ is a 2-cluster.)
Then the spectral shift function $\sigma$ is $\Cinf$ in $\Real\setminus(\Lambda
\cup\pspec(H))$, and it is a classical
symbol at infinity (i.e.\ outside a compact
set) with a complete asymptotic expansion:
\begin{equation*}
\sigma(\lambda)\sim\lambda^{\frac{n}{2}-3}\sum_{j=0}^\infty c_j\lambda^{-j},
\ c_0=C_0(n)\sum_{a,b:a\neq b}\int_{\Rn} V_a V_b\,dg.
\end{equation*}
\end{thm}

Note that $\sigma$ decays one order faster than in 2-body scattering, and
two orders faster than Weyl's law on compact manifolds. This is
because $\phi(H_0)+\sum_{\#b=2}(\phi(H_b)-\phi(H_0))$ is, in a high-energy
sense, closer to $\phi(H)$ than $\phi(H_0)$ is
to $\phi(H)$ in two body scattering.
If not all interactions are pair interactions, the
order of the leading term changes, namely becomes $\lambda^{\frac{n}{2}-2}$
as in 2-body scattering.

The proof of this theorem relies on the propagation of singularities, applied
to the Schwartz kernel of the resolvent, $R(\lambda+i0)$.
(In fact, the theorem should generalize to symbolic potentials, but the
proof would require a more precise microlocalization than provided by $\WFSc$.)
So we now
turn to the positive commutator methods that prove this.

\section{Microlocal positive commutator estimates}\label{sec:positive}
First I sketch, somewhat vaguely, the idea of positive commutator estimates.
So suppose that we want to obtain estimates on the solutions of $Pu=f$,
where $f$ is known, and is `nice', and $P$ is self-adjoint. Suppose that
we can construct an operator $A$ which is self-adjoint and is such that
\begin{equation}\label{eq:pos-comm-est-0}
i[A,P]=B^*B+E.
\end{equation}
Here $B^*B$ is the positive term, giving the name to the estimate.
The point is that we can estimate $Bu$ in terms of $Eu$.
Indeed, at least formally,
\begin{equation*}
\langle u,i[A,P]u\rangle=\langle u,B^*Bu\rangle+\langle u,Eu\rangle
=\|Bu\|^2+\langle u,Eu\rangle.
\end{equation*}
On the other hand,
\begin{equation*}
\langle u,i[A,P]u\rangle=\langle u,iAPu\rangle-\langle u,iPAu\rangle
=\langle u,iAPu\rangle+\langle iAPu,u\rangle=2\re\langle u,iAPu\rangle.
\end{equation*}
Combining these yields
\begin{equation}\label{eq:pos-comm-est}
\|Bu\|^2\leq 2|\re\langle u,iAPu\rangle|+|\langle u,Eu\rangle|.
\end{equation}
This means that $Bu$ can be estimated in terms of $Pu$, which is known
from the PDE, and $Eu$, on which we need to make assumptions. The typical
application is that $E$ is supported in one region of phase space and $B$
in another, in which case we can {\em propagate} estimates of $u$.

In fact, one can also apply this estimate if one does not know a priori that
$Bu\in L^2$. Namely, an approximation argument gives that if $Pu$ and $Eu$
are in appropriate spaces so that the right hand side of
\eqref{eq:pos-comm-est} makes sense, then $Bu\in L^2$, and
\eqref{eq:pos-comm-est} holds. Considering pseudo-differential operators
$A$ of various orders, this means that we obtain microlocal weighted Sobolev
estimates for $u$.
Also, typically one has an error term $F$, i.e.\ $i[A,P]=B^*B+E+F$, but
$F$ is `lower order' in some sense. Thus, $|\langle u,Fu\rangle|$
is added to the right hand side of \eqref{eq:pos-comm-est}, but
being `lower order' means that $|\langle u,Fu\rangle|$ automatically makes
sense, hence is irrelevant when proving that $Bu\in L^2$.

In fact, this method also yields estimates for the resolvent very directly.
Since for $t\in\Real$, $i[A,P-it]=i[A,P]$, and
\begin{equation*}\begin{split}
&\langle u,i[A,P-it]u\rangle=\langle u,iA(P-it)u\rangle-\langle u,i(P-it)Au\rangle\\
&\hspace{1cm}=\langle u,iA(P-it)u\rangle+\langle iA(P+it)u,u\rangle
=2\re\langle u,iA(P-it)u\rangle-2t\langle Au,u\rangle.
\end{split}\end{equation*}
Thus, we deduce
\begin{equation}\label{eq:pos-comm-est-t}
\|Bu\|^2+2t\langle Au,u\rangle\leq 2|\re\langle u,iAPu\rangle|+|\langle u,Eu\rangle|.
\end{equation}
This is in particular an estimate for $\|Bu\|$ provided that $t\geq 0$ and
$A$ is positive. Here we may take $u=u_t=(P-it)^{-1}f$, defined for
$t>0$, say, and we find a uniform estimate for $Bu_t$ as $t\to 0$.

The question is thus how one can produce operators $A$ which have
a positive commutator with $P$ as above. First, we recall how this happens
in the scattering calculus. Namely, if $A\in\Psisc^{m,l}(\Xb)$,
$P\in\Psisc^{m',l'}(\Xb)$ then $[A,P]\in\Psisc^{m+m'-1,l+l'+1}(\Xb)$
and
\begin{equation*}
\sigma_{m+m'-1,l+l'+1}(i[A,P])=H_a p=-H_p a,\ a=\sigma_{m,l}(A),
\ p=\sigma_{m',l'}(P),
\end{equation*}
where $H_a$ is the Hamilton vector field of $a$, $H_p$ the Hamilton vector
field of $p$. So modulo lower terms, which I ignore here and which are
easy to deal with, we need to arrange that
\begin{equation}\label{eq:Ham-vf-pos}
H_p a=-b^2+e,
\end{equation}
and then take $B,E$ with $\sigma(B)=b$, $\sigma(E)=e$. Indeed, under these
assumptions \eqref{eq:pos-comm-est} shows that $\|Bu\|$ can be estimated in
terms of $Pu$ and $Eu$. That is, $u$ microlocally
on $\supp b$ is estimated by $u$
on $\supp e$ (and $Pu$) in the precise sense described in the
next paragraph, so we can {\em propagate} estimates of $u$ from $\supp e$
to $\supp b$. (Incidentally, this is a good example of the $F$ term: only
the principal symbols of $E$ and
$B$ are specified. Take any $E$ and $B$ with these principal symbols,
$F=i[A,P]-B^*B-E$ is lower order.)

This can be used in a very straightforward manner to obtain bounds on
$\WFsc(u)$. Namely, one works with `relative wave front sets', relative
to $x^sH^r=H^{r,s}$, that is. Thus, for $\Xb=\overline{\Real^n}$,
$(z,\zeta)\nin\WFsc^{r,s}(u)$ means that
there is a cutoff function $\phi\in\Cinf_c(\Real^n)$ with $\phi(\zeta)\neq 0$
such that $\Fr^{-1}\phi\Fr u$ is in $H^{r,s}$ in an open cone around $z$.
But this is equivalent to the existence of some $Q\in\Psisc^{r,-s}(\Xb)$
such that $\sigma(Q)(z,\zeta)\neq 0$ and $Qu\in L^2$. Note that
$\sigma(Q)(z,\zeta)\neq 0$ means that $Q$ is elliptic at $(z,\zeta)$.
So if we find $A\in\Psisc^{m,l}(\Xb)$, and consequently
$B\in\Psisc^{\frac{m-1}{2},\frac{l+1}{2}}(\Xb)$,
$E\in\Psisc^{m-1,l+1}(\Xb)$ as above, then the conclusion is that (if
$Pu\in\dCinf(\Xb)$)
\begin{equation*}
\WFsc^{\frac{m-1}{2},-\frac{l+1}{2}}(u)\cap\supp e=\emptyset
\Rightarrow\WFsc^{\frac{m-1}{2},-\frac{l+1}{2}}(u)\cap\supp b=\emptyset.
\end{equation*}
In scattering theory $m$ is usually irrelevant by standard elliptic regularity.
Thus, one iteratively reduces $l$, proving that $\supp b$ is disjoint
from the wave front set with respect to more and more decaying Sobolev spaces.
(In fact, $b$ is shrunk slightly during the iterative procedure for
technical reasons.)

I now illustrate how to prove the propagation of singularities at $\pa\Xb$
for real principal type $P\in\Psisc^{m,0}(\Xb)$. For example, we may take
$P=H-\lambda$, $\lambda>0$, and $H$ is a two-body Hamiltonian.
(Note that microlocal elliptic regularity is
the consequence of the standard microlocal parametrix construction.)
We thus want to prove that if $Pu\in\dCinf(\Xb)$ (or a microlocal version of
it holds), $\bar\xi\in\Tsc^*_{\pa\Xb}\Xb$
and there is a point on the backward bicharacteristic through $\bar\xi$
which is not in $\WFsc(u)$, then $\bar\xi\nin\WFsc(u)$. In fact,
by a simple argument it suffices to prove that there exists a neighborhood $U$
of $\bar\xi$ such that if there is a point $\tilde\xi$ in $U$ which is also
on the backward bicharacteristic through $\bar\xi$ and
which is not in $\WFsc(u)$, then $\bar\xi\nin\WFsc(u)$.

The standard proof proceeds via linearization of $H_p$, see
\cite{Hormander:Existence}. Thus, first
note that $x^{-1}H_p$ is a smooth vector field on $\Tsc^*\Xb$ which
is tangent to the boundary. (For example, for Euclidean two-body Hamiltonians,
\begin{equation*}
x^{-1}H_p=2|z|\zeta\cdot\pa_z=2\sum_j\zeta_j |z|\pa_{z_j},
\end{equation*}
and $|z|\pa_{z_j}$ is a smooth vector field tangent to $\pa\Xb$,
i.e.\ it is in $\Vb(\Xb)$.) Thus, given any point
$\bar\xi\in\Tsc^*_{\pa\Xb}\Xb$
one can introduce local
coordinates $(q_1,\ldots,q_{2n-1})=(q_1,q'')$ on $\Tsc^*_{\pa\Xb}\Xb$ centered
at $\bar\xi$ such that
$x^{-1}H_p=\pa_{q_1}$ at $\Tsc^*_{\pa\Xb}\Xb$. Thus, bicharacteristics
at $\pa\Xb$ are curves $q''=\text{constant}$.
Now let $\chi_1\in\Cinf_c(\Real_{q_1})$ and
$\chi_2\in\Cinf_c(\Real^{2n-2}_{q''})$
be smooth functions supported near $0$ with the property that
\begin{equation*}
\chi_1'=-b_1^2+e_1,
\end{equation*}
$b_1,e_1\in\Cinf_c(\Real)$, and $\supp e_1\subset (-\infty,0)$.
Let
\begin{equation*}
a=\chi_1\chi_2^2,\ b=b_1\chi_2,\ e=e_1\chi_2^2.
\end{equation*}
Then \eqref{eq:Ham-vf-pos} holds. In fact, we can even allow weights
and take
\begin{equation*}
a_s=x^s\chi_1\chi_2^2, s\in\Real,
\end{equation*}
since $(x^{-1}H_p x^s)\chi_1$ can be absorbed in $x^s(x^{-1}H_p\chi_1)$
by choosing $\chi_1'$
large compared to $\chi_1$. This gives microlocal weighted
$L^2$ estimates in $x^{-s-1/2}L^2$.
The iterative argument, in which we gradually let $s\to-\infty$,
then allows one to conclude that
\begin{equation*}
\supp e\cap\WFsc(u)=\emptyset\Mand\supp a\cap\WFsc(Pu)=\emptyset\Rightarrow
\{b>0\}\cap\WFsc(u)=\emptyset.
\end{equation*}
By choosing $\supp\chi_1$ and $\supp\chi_2$ appropriately, we may arrange
that $e$ is supported near $\tilde\xi$ so that
$\supp e\cap\WFsc(u)=\emptyset$, and so that $b(\bar\xi)>0$, as shown
below.

\begin{figure}
\begin{center}
\mbox{\epsfig{file=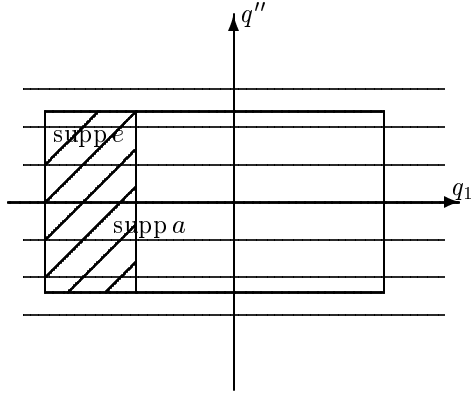}}
\end{center}
\caption{$\supp a$ superimposed on the linearized Hamilton flow.
$\supp e$ is the shaded region on the left.}
\label{fig:rptlin}
\end{figure}

There are several directions from here. One can use finer notion of
regularity, such as Lagrangian regularity, which would correspond to using
$\chi_2$ that vanishes simply on a Lagrangian submanifold, or such as
regularity at radial points, which is the subject of a joint paper
with Andrew Hassell and Richard Melrose \cite{Hassell-Melrose-Vasy:Spectral}.

Here I will talk about a rougher version, namely what happens if the
bicharacteristic `flow' is more complicated, e.g.\ in the presence of
boundaries or corners for the wave equation
\cite{Melrose-Sjostrand:I, Lebeau:Propagation} or many-body scattering.
In fact, here I will not explain the detailed behavior of bicharacteristics,
rather just show how to microlocalize positive commutator estimates in
a versatile fashion. This method goes back to the work of
Melrose and Sj\"ostrand \cite {Melrose-Sjostrand:I}.

The main point is that if we cannot put the operator $P$, or
at least its Hamilton vector field $H_p$ in a model form, the previous
construction will not work. Indeed, unless $H_p\chi_2=0$, $H_p(\chi_1\chi_2)$
will always yield a term $\chi_1 H_p\chi_2$, which cannot be controlled
by $(H_p\chi_1)\chi_2$: the problem being near the boundary of $\supp\chi_2$.
So instead use a different form of localization.
First let $\eta\in\Cinf(\Tsc^*\Xb)$ be a function with
\begin{equation*}
\eta(\bar\xi)=0,\ H_p\eta(\bar\xi)>0.
\end{equation*}
Thus, $\eta$ measures propagation along bicharacteristics, e.g.\ $\eta=q_1$
in the above example would work, but so would many other choices.
We will use a function $\omega$ to localize near putative bicharacteristics.
This statement is deliberately vague; at first
we only assume that $\omega\in\Cinf(\Tsc^*\Xb)$ is the sum of the squares of
$\Cinf$ functions $\sigma_j$, $j=1,\ldots,l$,
with non-zero differentials at $\bar\xi$ such that $d\eta$ and
$d\sigma_j$, $j=1,\ldots,l$, span $T_{\bar\xi}\Tsc^*_{\pa\Xb}\Xb$.
Such a function $\omega$ is non-negative
and it vanishes
quadratically at $\bar\xi$, i.e.\ $\omega(\bar\xi)=0$ and
$d\omega(\bar\xi)=0$. An example is $\omega=q_2^2+\ldots+q_{2n-1}^2$
with the notation from before, but again there are many other possible choices.
We now consider a family symbols, parameterized by constants 
$\delta\in(0,1)$, $\ep\in(0,\delta]$, of the form
\begin{equation*}
a=\chi_0(2-\frac{\phi}{\ep})\chi_1(\frac{\eta+\delta}{\ep\delta}+1),
\end{equation*}
where
\begin{equation*}
\phi=\eta+\frac{1}{\ep}\omega,
\end{equation*}
$\chi_0(t)=0$ if $t\leq 0$, $\chi_0(t)=e^{-1/t}$ if $t>0$,
$\chi_1\in\Cinf(\Real)$, $\supp\chi_1\subset[0,+\infty)$,
$\supp\chi_1'\subset[0,1]$.
Although we do not do it explicitly here, weights such as $x^s$ can
be accommodated for any $s\in\Real$, by replacing the factor
$\chi_0(2-\frac{\phi}{\ep})$ by $\chi_0(A_0^{-1}(2-\frac{\phi}{\ep}))$
and taking $A_0>0$ large.

We analyze the properties of $a$ step by step. First, note
that $\phi(\bar\xi)=0$, $H_p\phi(\bar\xi)=H_p\eta(\bar\xi)>0$,
and $\chi_1(\frac{\eta+\delta}{\ep\delta}+1)$ is identically $1$
near $\bar\xi$, so $H_p a(\bar\xi)<0$. Thus, $H_p a$ has
the correct sign, and is in particular non-zero, at $\bar\xi$.

Next,
\begin{equation*}
\xi\in\supp a\Rightarrow\phi(\xi)\leq2\ep\Mand
\eta(\xi)\geq-\delta-\ep\delta.
\end{equation*}
Since $\ep<1$, we deduce that in fact $\eta=\eta(\xi)\geq -2\delta$.
But $\omega\geq 0$, so $\phi=\phi(\xi)\leq 2\ep$ implies that
$\eta=\phi-\ep^{-1}\omega\leq \phi\leq 2\ep$. Hence, $\omega=\omega(\xi)
=\ep(\phi-\eta)\leq 4\ep\delta$. Since $\omega$ vanishes quadratically
at $\bar\xi$, it is useful to rewrite the estimate as $\omega^{1/2}\leq
2(\ep\delta)^{1/2}$. Combining these, we have seen that on $\supp a$,
\begin{equation}\label{eq:supp-a-est}
-\delta-\ep\delta\leq \eta\leq 2\ep\Mand \omega^{1/2}\leq 2(\ep\delta)^{1/2}.
\end{equation}
Moreover, on $\supp a\cap\supp\chi_1'$, 
\begin{equation*}
-\delta-\ep\delta\leq \eta\leq -\delta
\Mand \omega^{1/2}\leq 2(\ep\delta)^{1/2}.
\end{equation*}
Note that given any neighborhood $U$
of $\bar\xi$, we can thus make $a$ supported in $U$ by
choosing $\ep$ and $\delta$ sufficiently small.
Below we illustrate the parabola shaped region given by $\supp a$.

\begin{figure}
\begin{center}
\mbox{\epsfig{file=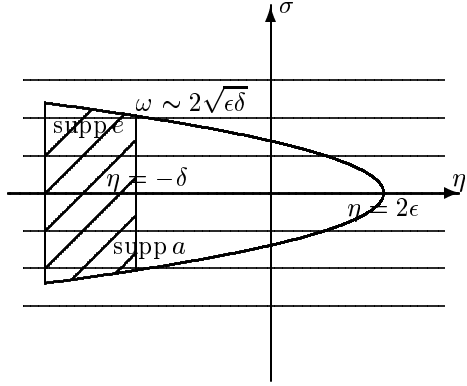}}
\end{center}
\caption{$\supp a$ in $(\eta,\sigma)$ coordinates.
$\supp e$ is again the shaded region on the left.}
\label{fig:rptpar2}
\end{figure}

Note that as $\ep\to 0$, but $\delta$ fixed,
the parabola becomes very sharply localized at $\omega=0$. In particular,
for very small $\ep>0$ we obtain a picture quite analogous to letting
$\supp\chi_2\to\{0\}$ in Figure~\ref{fig:rptlin}.

So we have shown that $a$ is supported near $\bar\xi$. We define
\begin{equation*}
e=\chi_0(2-\frac{\phi}{\ep})H_p\chi_1(\frac{\eta+\delta}{\ep\delta}+1),
\end{equation*}
so the crucial question is whether $H_p\phi\geq 0$ on $\supp a$. Note that
choosing $\delta_0\in(0,1)$ and $\ep_0\in(0,\delta_0)$ sufficiently small,
one has
$H_p\eta\geq c_0>0$ where $|\eta|\leq 2\delta_0$,
$\omega^{1/2}\leq 2(\ep_0\delta_0)^{1/2}$. So $H_p\phi\geq 0$ on $\supp a$,
provided that $|H_p\omega|\leq \frac{c_0}{2}\ep$ there.

But being a sum of squares of functions with non-zero differentials,
$H_p\omega$ vanishes at $\omega=0$ and satisfies
$|H_p\omega|\leq C\omega^{1/2}$. Due to \eqref{eq:supp-a-est},
we deduce that $|H_p\omega|\leq 2C(\ep\delta)^{1/2}$. So
$|H_p\omega|\leq \frac{c_0}{2}\ep$ holds if $\frac{c_0}{2}\ep\geq
2C(\ep\delta)^{1/2}$, i.e.\ if $\ep\geq C'\delta$ for some constant $C'>0$
independent of $\ep$, $\delta$. Note that this constraint on $\ep$,
i.e.\ that it cannot be too small, gives very rough localization:
the width of the parabola at $\eta=-\delta$ is roughly
$\omega^{1/2}\sim\delta$, i.e.\ it is very wide, and in particular insufficient
to prove the propagation of singularities along the bicharacteristics.
The reason is simple:
our localizing function $\omega$ has no relation to $H_p$, so we cannot
expect a more precise estimate. Note, however, that the estimate is still
non-trivial! Indeed, it shows that singularities propagate in the
sense that $\bar\xi$ cannot be an isolated point of $\WFsc(u)$.
(We required $\ep\in(0,\delta]$ beforehand, but in fact we could have
dealt with $\ep\leq\mu\delta$, even if $\mu>1$, if we localized slightly
differently.)

We need to adapt $\omega$ to $H_p$ to get a better estimate. If
we linearize $H_p$ as above, and take $\omega=q_2^2+\ldots+q_{2n-1}^2$,
then $H_p\omega=0$ and any $\ep>0$ works. Thus, in this case, we can
prove propagation of singularities much like by the previous, simpler,
construction.

However, we do not need such a strong relationship to $H_p$. Suppose instead
that we merely get $\omega$ `right' at $\bar\xi$, in the sense that
$\omega=\sum \sigma_j^2$ and $H_p\sigma_j(\bar\xi)=0$. Then
$|H_p\sigma_j|\leq C_0(\omega^{1/2}+|\eta|)$, so
$|H_p\omega|\leq C \omega^{1/2}(\omega^{1/2}+|\eta|)$. Using
\eqref{eq:supp-a-est}, we deduce that $|H_p\omega|\leq \frac{c_0}{2}\ep$
provided that $\frac{c_0}{2}\ep\geq C''(\ep\delta)^{1/2}\delta$, i.e.\ that
$\ep\geq C'\delta^3$ for some constant $C'$ independent of $\ep$, $\delta$.
Now the size of the parabola at $\eta=-\delta$ is roughly
$\omega^{1/2}\sim\delta^2$, i.e.\ we have localized along a single
direction, namely the direction of $H_p$ at $\bar\xi$. By a relatively
simple argument, one can piece together such estimates (i.e.\ where
the direction is correct `to first order') and deduce the
propagation of singularities. We emphasize that the lower bound for $\ep$
is natural. Indeed, with $q_j$ as above, we may take $\sigma_j$ e.g. to be
$\sigma_j=q_j+q_1^2$, $j\geq 2$. The bicharacteristics are
$q_j=\text{constant}$, but we are localizing near $\sigma_j=\text{constant}$,
and at $\eta=-\delta$ these differ by $\delta^2$. So any localization better
than $\omega^{1/2}\sim\delta^2$ would in fact contradict the propagation
of singularities!

\begin{figure}
\begin{center}
\mbox{\epsfig{file=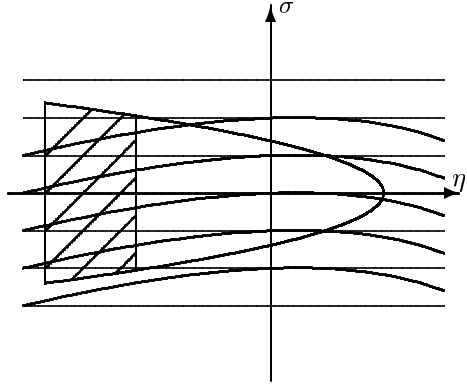}}
\end{center}
\caption{Bicharacteristics and $\supp a$. The labels from
Figure~\ref{fig:rptpar2} have been removed to make the picture
less cluttered. The straight horizontal lines are the $\sigma=\text{constant}$
lines, while the nearby parabolae are the bicharacteristics.}
\label{fig:rptpar3}
\end{figure}

The microlocal positive commutator estimates in many-body scattering arise
by this method. In particular, one can take $\eta=\frac{z\cdot\zeta}{|z|}$,
which is the radial component of the momentum. The function $\omega$ needs
to be $\pi$-invariant, so if $\bar\xi\in\Tsc^*\Xb_a$, it involves
functions on $\Tsc^*\Xb_a$ as well as $Z^a=\frac{z^a}{|z|}$ and $\eta$.
The only additional argument needed is to show that the commutator is
indeed positive,
which has to be understood in an operator sense. Thus, the key point
is that the commutator of $H-\lambda$ with a quantization $B$ of $\eta$
is positive, modulo compact terms, when localized at $\lambda$ in the
spectrum of $H$ and microlocalized away from the radial set $\calR(\lambda)$.
Note that, as usual, there is nothing to prove at $\calR(\lambda)$, since
each point in it is the image of a maximally extended generalized broken
bicharacteristic.

This positivity can be proved directly by showing
that the indicial operators of the commutator are positive away from
$\calR(\lambda)$, which follows from an iterative argument.
However, it also reduces to the Mourre estimate, involving the generator
of dilations $A=\frac{1}{2}(\cdot D_z+D_z\cdot z)$, which has
principal symbol at $\zeta\cdot z$.
The Mourre estimate states the following. Suppose that $\lambda\nin\Lambda$ and
$\ep>0$. Then there is a $\delta>0$ such that if
$\phi\in\Cinf_c(\Real)$ is supported in $(\lambda-\delta,\lambda+\delta)$
then there exists $K'\in\PsiSc^{-\infty,1}(\Xb,\calC)$ such that
\begin{equation}\label{eq:Mourre-0}
\phi(H)i[A,H]\phi(H))\geq 2(d(\lambda)-\ep)\phi(H)^2+K',
\end{equation}
where
\begin{equation*}
d(\lambda)=\inf\{\lambda-\lambda':\lambda'\leq\lambda,\ \lambda'\in\Lambda\}
\geq 0
\end{equation*}
is the distance of $\lambda$ to the next threshold below it if $\lambda\geq
\inf\Lambda$, and $d(\lambda)$ an arbitrary positive number if
$\lambda<d(\lambda)$. Since $d(\lambda)>0$ if $\lambda\nin\Lambda$,
\eqref{eq:Mourre-0} is indeed a positive commutator estimate, which
does not even have a `negative' term $E$, unlike \eqref{eq:pos-comm-est-0}.
The Mourre estimate, originating in \cite{Mourre-Absence}, has been well
understood since the work of Perry, Sigal and Simon
\cite{Perry-Sigal-Simon:Spectral} and Froese and Herbst~\cite{FroMourre}.
Here I just outline the argument in the simplest case, namely if no
proper subsystem has any $L^2$-eigenvalues.

In this simplest case, the argument of \cite{FroMourre} (see also
\cite{Derezinski-Gerard:Scattering} and \cite{Vasy:Exponential})
proceeds as follows.
In order to prove \eqref{eq:Mourre-0}, one only
needs to show that for all $b$, the corresponding indicial operators
satisfy the corresponding inequality, i.e.\ that
\begin{equation}\label{eq:Mourre-1}
\phi(\hat H_b)i\widehat{[A,H]}_b\phi(\hat H_b)
\geq 2(d(\lambda)-\ep)\phi(\hat H_b)^2.
\end{equation}
(This means that the operators on the two sides, which are families of
operators
on $X^b$, depending on $(y_b,\zeta_b)\in\sct_{C_b}\Xb$, satisfy the inequality
for all $(y_b,\zeta_b)\in\sct_{C_b}\Xb$.)
It is convenient to assume that $\phi$ is identically $1$ near
$\lambda$; if \eqref{eq:Mourre-1} holds for such $\phi$, it holds
for any $\phi_0$ with slightly smaller support, as follows by multiplication
by $\phi_0(\hat H_b)$ from the left and right.

Note that for $b=0$ the estimate
certainly holds: it comes from the Poisson bracket formula
in the scattering calculus, or from a direct computation yielding
$i\widehat{[A,H]}_0=2\Delta_{g_0}$. Hence, if the the localizing factor
$\phi(\hat H_0)=\phi(|\zeta|^2)$ is supported in
$(\lambda-\delta,\lambda+\delta)$ and $\lambda>0$,
then \eqref{eq:Mourre-1} holds even
with $d(\lambda)-\ep$ replaced with $\lambda-\delta$. Note that
$\lambda\geq d(\lambda)$, if $\lambda>0$, since $0$ is a threshold of $H$.
On the other hand, if $\lambda<0$, both sides of \eqref{eq:Mourre-1}
vanish for $\phi$ supported near $\lambda$, so the inequality holds
trivially.

In general, we may assume inductively that at all clusters $c$ with
$C_c\subsetneq C_b$, i.e.\ $X^b\subsetneq X^c$, \eqref{eq:Mourre-1} has
been proved with $\phi$ replaced by a cutoff $\phit$ and $\ep$ replaced
by $\ep'$, i.e.\ we may assume that for all $\ep'>0$ there exists $\delta'>0$
such that for all $c$ with $C_c\subsetneq C_b$, and for all
$\phit\in\Cinf_c(\Real;[0,1])$ supported in
$(\lambda-\delta',\lambda+\delta')$,
\begin{equation}\label{eq:Mourre-1-0}
\phit(\hat H_c)i\widehat{[A,H]}_c\phit(\hat H_c)
\geq 2(d(\lambda)-\ep')\phit(\hat H_c)^2.
\end{equation}
But these are exactly the indicial operators of
$\phit(\hat H_b)i\widehat{[A,H]}_b\phit(\hat H_b)$, so, as
discussed in \cite[Proposition~8.2]{Vasy:Propagation-Many},
\eqref{eq:Mourre-1} implies that
\begin{equation}\label{eq:Mourre-8}
\phit(\hat H_b)i\widehat{[A,H]}_b\phit(\hat H_b)
\geq 2(d(\lambda)-\ep')\phit(\hat H_b)^2+K_b,
\ K_b\in\PsiSc^{-\infty,1}(X^b,\calC^b).
\end{equation}
This implication relies on a square root construction in
the many-body calculus, which is particularly simple in this case.

Now, we first multiply \eqref{eq:Mourre-8}
through by $\phi(H)$ from both the left and the
right. Recall that we use coordinates $(z_b,z^b)$ on $X_b\oplus X^b$ and
$(\zeta_b,\zeta^b)$ are the dual coordinates.
We remark that $\hat H_b=|\zeta_b|^2+H^b$, so
if $\lambda-|\zeta_b|^2$ is not an eigenvalue of $H^b$,
then as $\supp\phi\to\{\lambda\}$,
$\phi(H^b+|\zeta_b|^2)\to 0$ strongly, so as $K_b$ is compact,
$\phi(H^b+|\zeta_b|^2)K_b\to 0$ in norm; in particular it
can be made to have norm smaller than $\ep'-\ep>0$. After multiplication
from both sides by $\phi_1(\hat H_b)$, with $\phi_1$ having even smaller
support, \eqref{eq:Mourre-1}
follows (with $\phi_1$ in place of $\phi$),
with the size of $\supp\phi_1$ a priori depending on $\zeta_b$.
However, $i\phi_1(\hat H_b)\widehat{[A,H]}_b\phi_1(\hat H_b)$
is continuous in $\zeta_b$ with values in bounded operators on $L^2(X^b)$,
so if \eqref{eq:Mourre-1} holds at one value
of $\zeta_b$, then it holds nearby. Moreover, for large $|\zeta_b|$ both sides
vanish as $\hat H_b=H^b+|\zeta_b|^2$, with $H^b$ bounded below, so
the estimate is in fact uniform
if we slightly increase $\ep>0$.

In general, the proof requires to treat the range of $E$, the spectral
projection of $H^b$ to $\{\lambda\}$, separately. Roughly, the
positivity estimate on the range of $E$ comes from the virial theorem,
$iE[z^bD_{z^b},H^b]E=0$, which is formally clear, and is easy to prove.
Thus,
\begin{equation*}
iE[A,H_b]E=iE[z^bD_{z^b},H^b]E+iE[z_bD_{z_b},\Delta_{X_b}]E
=i[z_bD_{z_b},\Delta_{X_b}]E,
\end{equation*}
and the commutator $i\phi(H_b)
[z_bD_{z_b},\Delta_{X_b}]\phi(H_b)$ is easily computed to be
positive. Of course, there are also cross-terms that need to be
considered, but they can be estimated by Cauchy-Schwartz estimates,
see \cite{FroMourre} or \cite{Vasy:Exponential}.

I refer to \cite{Vasy:Bound-States}
and \cite{Vasy:Propagation-Many} for the detailed arguments proving
propagation of singularities in the many-body setting,
and to \cite[Appendix]{Vasy-Wang:Smoothness} for
weaker estimates with simplified proofs.

\section{Asymptotic completeness}\label{sec:ac}
Asymptotic completeness (AC)
is an $L^2$-based statement describing the long-term
behavior of solutions of the Schr\"odinger equation. In the short-range
setting it says that
for any $\phi\in L^2(X_0)$ in the range of $\Id-E_{\mathrm{pp}}$,
$E_{\mathrm{pp}}$ being the projection onto the bound states of $H$
(i.e.\ onto its $L^2$-eigenfunctions),
there exist $\phi_\alpha\in L^2(X_a)$ such
that
\begin{equation*}
\|e^{-iHt}\phi-\sum_\alpha e^{-iH_a t}(\phi_\alpha\otimes\psi_\alpha)\|\to 0
\ \text{as}\ t\to+\infty.
\end{equation*}
In the long-range setting, $e^{-iH_a t}$ must be somewhat modified.
After the ground-breaking work of Enss~\cite{Enss:Asymptotic, Enss:Quantum}, AC
was first proved by Sigal and Soffer \cite{Sigal-Soffer:N} in
the short-range setting (see Graf's paper \cite{Graf:Asymptotic} for
a different proof), and later by Derezi\'nski
\cite{Derezinski:Asymptotic}, and also by Sigal and Soffer
\cite{Sigal-Soffer:Long-range, Sigal-Soffer:Asymptotic},
in the long-range setting. In the short range case the main ingredient
is equivalent to certain estimates of the resolvent at the radial sets
in a sense that I now describe. In the long-range setting, as already
in two-body scattering, additional constructive steps are needed.
The estimates, in a different language, appeared first in the
work of Sigal and Soffer \cite{Sigal-Soffer:N}. I hope that the
following discussion makes it clearer how they relate to the propagation of
singularities.

While asymptotic completeness gives a complete long-term
$L^2$-description of solutions
of the Schr\"odinger equation, the question remains whether an analogous
description exists on other spaces, such as weighted $L^2$-spaces. For
example, if $\phi$ is Schwartz, are the $\phi_\alpha$ Schwartz?
Or dually, starting with a tempered distribution $\phi$, are there
tempered distributions $\alpha$ such that the convergence holds, as
$t\to+\infty$, in a suitable sense? A different point of view
is the parameterization of generalized eigenfunctions of $H$ using
the Poisson operators $P_{\alpha,+}(\lambda)$, and the analogues of
these questions can be asked there as well. The answer is affirmative in
the two-body setting (even in the geometric setting, see
\cite{RBMSpec, RBMZw}). However, as indicated by
the related issue of the mapping properties
of the scattering matrices, discussed at the end of this section,
it is unlikely that the same holds in the many-body setting. One can then
ask weaker question, e.g.\ whether it holds in weighted spaces $L^2_s$,
$s$ near $0$. Or, one may ask
whether one can give a precise description of the map
$\phi\mapsto\phi_\alpha$ e.g.\ as some sort of Fourier integral operator.

As a starting point of relating the propagation of singularities to AC,
we note that the propagation of singularities is
proved by showing its `relative' versions, i.e.\ that for any $l$,
$\WFSc^{*,l}(u)$
is also a union of maximally extended generalized broken bicharacteristics.
When considering the resolvent, first recall that for
$f\in\dCinf(\Xb)$, $R(\lambda+ i0)f\in H^{\infty,l}$ for all $l<-1/2$,
so we only need to find $\WFSc^{*,l}(R(\lambda+ i0)f)$ for $l\geq -1/2$.
Theorem~\ref{thm:lim-abs-WF} is also valid for $\WFSc^{*,l}$, i.e.\ the
following holds.

\begin{thm}\label{thm:lim-abs-WF-l}
If $\lambda\nin\Lambda$, then for $f\in\Sch(\Rn)$, $l\geq -1/2$,
$\WFSc(R(\lambda+ i0)f)$ is a subset of the image of $\calR_+(\lambda)$
under the {\em forward} generalized broken bicharacteristic relation.
\end{thm}

This result allows $u=R(\lambda+ i0)f$ {\em not} to lie in $H^{*,-1/2}$
on the image of $\calR_+(\lambda)$
under the forward generalized broken bicharacteristic relation.
This is a small set, but it is important to know whether $\WFSc^{*,l}(u)$
may indeed intersect the forward image of $\calR_+(\lambda)$. Of course,
we cannot expect an improvement at $\calR_+(\lambda)$, as shown already
by the example of the free Euclidean Laplacian. The crucial improvement
is the following estimate, due to Sigal and Soffer \cite{Sigal-Soffer:N}.

\begin{thm}\label{thm:lim-abs-WF-12}
If $\lambda\nin\Lambda$, then for $f\in\Sch(\Rn)$,
$\WFSc^{*,-1/2}(R(\lambda+ i0)f)\subset \calR_+(\lambda)$.
\end{thm}

\begin{rem*}
This theorem also has a time-dependent analogue. If
$u$ is a solution of the Schr\"odinger equation $(D_t+H)u=0$ with
$u|_{t=0}\in\Sch(X_0)$ then on the one hand
$u\in H^{\infty,l}(X_0\times\Real)$ for $l<-1/2$, on the other
hand $\WFSc^{*,-1/2}(u)\subset\calR$.
\end{rem*}

In fact, this theorem can be improved along the lines of the distributional
statement in Theorem~\ref{thm:lim-abs-WF}:

\begin{cor}\label{cor:lim-abs-WF-12}
Suppose that $\lambda\nin\Lambda$, $f\in H^{*,1/2}$ and
$\WFSc^{*,1/2+\ep}(f)\cap \calR_-(\lambda)
=\emptyset$ for some $\ep>0$.
Then $R(\lambda+i0)f=\lim_{t\to 0}R(\lambda+it)f$ exists
in $H^{*,-1/2-\ep'}(\Xb)$, $\ep'>0$,
and $\WFSc^{*,-1/2}(R(\lambda+i0)f)\subset
\calR_+(\lambda)$.
\end{cor}

Theorem~\ref{thm:lim-abs-WF-12}
can be proved rather simply. The main issue is how to obtain
a positive commutator at the radial point. Away from the radial
sets arbitrary weights can be accommodated by suitable construction,
as pointed out in the previous section. At radial points only
the weights can give positive commutators. Now, one has to use weights
$x^{-2l-1}$ to obtain estimates for $\WFSc^{*,l}$, and these
weights will give a commutator whose sign depends on that of $-2l-1$, hence
on whether $l>-1/2$, $l<-1/2$ or $l=-1/2$. It turns out that the sign is
correct for \eqref{eq:pos-comm-est-t} to be of use if $l<-1/2$; this
yields the limiting absorption principle. The sign is wrong if $l>-1/2$,
so no results can be expected then. In the borderline case $l=-1/2$, the
weight vanishes. The way to obtain a positive commutator is thus to consider
operators $A$ which are microlocally (a multiple of) the identity
near $\calR_+(\lambda)$. The commutator then vanishes microlocally near
$\calR_+(\lambda)$, which is reasonable since no estimate on $\WFSc^{*,-1/2}(u)$
can be expected there.

It is then straightfoward to construct $A$
so that \eqref{eq:pos-comm-est-t} can be used to prove
Theorem~\ref{thm:lim-abs-WF-12}. Indeed, it suffices to show
that on $\WFSc^{*,-1/2}(u)$, $\eta=\frac{z\cdot\zeta}{|z|}$
must satisfy $\lambda-\eta^2\in\Lambda$, for then the full statement
of the theorem follows by the propagation of singularities for
$\WFSc^{*,-1/2}(u)$. So we
proceed to prove
this simpler result, namely that if $\lambda-\bar\eta^2\nin\Lambda$ then
for any point $\xi$, $\eta(\xi)=\bar\eta$ implies $\xi\nin
\WFSc^{*,-1/2}(u)$.

To do so, we let
$a=\chi(\eta)$ where $\chi\in\Cinf_c(\Real)$, $\chi
\geq 0$,
is chosen so that $\chi\equiv 1$ on $[0,\bar\eta-\delta]$ for some $\delta>0$,
$\chi'\leq 0$ on $(0,\infty)$, $\chi'(\bar\eta)<0$,
and $t\in\supp\chi'$ implies that $\lambda-t^2\nin\Lambda$. This can
be arranged as $\Lambda$ is closed. We can
further make sure that $\sqrt{-\chi'}$ is $\Cinf$ on $(0,\infty)$.
Then the positive commutator methods outlined show the commutator
of $A$, a quantization of $a$, with $H-\lambda$ is positive, in the
region $\eta>0$, yielding the estimate that proves the theorem.
We remark that partial microlocalization, using functions of $\eta$,
hes been used extensively in many-body scattering, especially
by G\'erard, Isozaki and Skibsted \cite{GerComm, GIS:N-body}
and Wang \cite{XPWang}, to obtain
partially microlocal statements such as radiation conditions and
uniqueness statements
\cite{IsoUniq, IsoRad},
and indeed to prove the smoothness of 2-cluster to 2-cluster
scattering matrices \cite{Skibsted:Smoothness}.

It turns out that there is an even simpler way of proving
Theorem~\ref{thm:lim-abs-WF-12}, or indeed a stronger statement, which is
due to Yafaev \cite{Yaf}. His estimate states
that in a neighborhood of $C_{a,\reg}$, where we
write $y_a$ for the coordinates $z_a/|z_a|$ along $C_{a,\reg}$,
$xD_{y_a} R(\lambda+i0)f$ is in $H^{*,-1/2}(\Xb)$. Since
the principal symbol of $xD_{y_a}$ is invertible on $(T^*X_{a,\reg}\cap
\dChar(\lambda))\setminus \calR(\lambda)$, this result implies
Theorem~\ref{thm:lim-abs-WF-12}. Yafaev's proof relies on a simple and
explicit commutator calculation, which allows one to deal with various
error terms that one may, a priori, expect. However, exactly because of its
explicit nature, it is presumably hard to generalize to more geometric
settings, while the argument we sketched does not face this difficulty.

As discussed by Yafaev
\cite{Yaf} in the usual time-dependent version, short-range
asymptotic clustering, hence asymptotic completeness,
are relatively easy consequences of Corollary~\ref{cor:lim-abs-WF-12},
and we refer to \cite{Yaf} and \cite{Derezinski-Gerard:Scattering}
for more details. However, it is worth pointing out that the reason
why Coulomb-type potentials (i.e.\ those in $S^{-1}$) are not `short-range'
is that the Hamilton vector field in some
subsystem vanishes at radial points. This degeneracy makes even the
subprincipal term important in describing the precise behavior
of generalized eigenfunctions microlocally near this point.

Before turning to scattering theory on symmetric spaces, we note the
implications of Theorem~\ref{thm:lim-abs-WF-12} for the scattering
matrices. Previously, $S_{\alpha\beta}(\lambda)$ was only defined
as a map $S_{\alpha\beta}(\lambda):\Cinf_c(S_{a,\reg})\to\dist(S_{b,\reg})$.
Indeed, part of the broken bicharacteristic relation connects
$\calR_+(\lambda)$ with its image, and this can a priori give
a singularity in the kernel of $S_{\alpha\beta}(\lambda)$ of the kind
that does not even allow one to conclude that
$S_{\alpha\beta}(\lambda):\Cinf_c(S_{a,\reg})\to\Cinf(S_{b,\reg})$.
The pairing formula, \eqref{eq:S-pairing}, combined with
Theorem~\ref{thm:lim-abs-WF-12}, show that in fact
\begin{equation}\label{eq:S-L2}
S_{\alpha\beta}(\lambda):L^2(S_a)\to L^2(S_b).
\end{equation}
It is an interesting question whether this can be improved if
we restrict $S_{\alpha\beta}(\lambda)$ to $\Cinf_c(S_{a,\reg})$. Namely,
except in special cases such as $N$-clusters and two-clusters, the
best known result is the trivial consequence of \eqref{eq:S-L2}:
\begin{equation*}
S_{\alpha\beta}(\lambda):\Cinf_c(S_{a,\reg})\to L^2(S_b).
\end{equation*}
(In the case of $N$-clusters and 2-clusters, the geometry of generalized
broken bicharacteristics gives
$S_{\alpha\beta}(\lambda):\Cinf_c(S_{a,\reg})\to\Cinf(S_{b,\reg})$.)
The putative improvement would have to be connected to an improvement
of Theorem~\ref{thm:lim-abs-WF-12}, namely to the
existence of {\em some} $l>-1/2$ such that
$\WFSc^{*,l}(R(\lambda+ i0)f)\subset \calR_+(\lambda)$. It would also
be connected to a better understanding of $R(\lambda\pm i0)$ at the
thresholds, in which direction Wang's paper \cite{Wang:Spectral} is
the only one I am aware of.

\section{Scattering on higher rank symmetric spaces}\label{sec:sl3}
In this section I discuss $\SL(N,\Real)/\SO(N,\Real)$, indeed, mostly
$\SL(3,\Real)/\SO(3,\Real)$.
The books \cite{Helgason:Groups}, \cite{Jost:Riemannian} and
\cite{Eberlein:Geometry} are good general references.
$N=2$ yields the hyperbolic plane $\mathbb{H}^2$, which is a rank one symmetric
space on which many aspects of analysis, such as the asymptotic behavior
of the resolvent kernel and the analytic continuation of the
resolvent are well understood. Indeed, these have been described
on asymptotically hyperbolic spaces by
Mazzeo and Melrose \cite{Mazzeo-Melrose:Meromorphic}
and Perry \cite{Perry:Laplace-I, Perry:Laplace-II}.

Higher rank symmetric spaces, such as $\SL(N)/\SO(N)$, $N\geq 3$,
are much less understood. For example, using results
of Harish-Chandra, and Trombi and Varadarajan (see \cite{Helgason:Groups}),
Anker and Ji
only recently obtained the leading order behavior of the Green's function
\cite{Anker:Forme, Anker-Ji:Comportement, Anker-Ji:Heat}.
Also, while spherical functions, which are most analogous to partial
plane-partial spherical waves in the Euclidean setting, have been analyzed by 
Harish-Chandra, Trombi and Varadarajan, and in particular their
analytic continuation is understood, the same cannot be said about
the Green's function. The analysis of spherical functions relies
on perturbation series expansions, much like in the proof of the
Cauchy-Kovalevskaya theorem, and it does not work well at the walls of
the Weyl chambers. Here I only illustrate some recent joint results
with Rafe Mazzeo \cite{Mazzeo-Vasy:Sl3, Mazzeo-Vasy:Sl3-Analytic},
that illuminate the connections with many-body
scattering, and in particular give rather direct results for the
resolvent.

First I describe the space $\SL(3)/\SO(3)$. The polar decomposition
states that any $C\in\SL(3)$ can be written uniquely as
$C=VR$, $V=(CC^t)^{1/2}$ is positive definite and has determinant $1$,
$R\in\SO(3)$. Thus, $\SL(3)/\SO(3)$ can be identified with the
set $M$ of positive definite matrices of determinant $1$;
this is a five-dimensional real analytic manifold.
The Killing form provides a Riemannian
metric $g$. The associated Laplacian $\Delta = \Delta_g$ gives a self-adjoint 
unbounded operator on $L^2(M,dg)$, with spectrum $[\lambda_0,+\infty)$, 
$\lambda_0=\frac{1}{3}$. Let $R(\lambda)=(\Delta-\lambda)^{-1}$ be the resolvent of 
$\Delta_g$, $\lambda\nin [\lambda_0,+\infty)$.

Fix a point $o \in M$, which
we may as well assume is the image of the identity matrix $I$ in
the identification above. The stabilizer subgroup $K_o$ (in the
natural $\SL(3)$ action on $M$) is isomorphic to $\SO(3)$.
The Green function $G_o(\lambda)$ with pole at $o$ and at eigenvalue $\lambda$ is, by 
definition $R(\lambda)\delta_o$. It is standard that $G_o$ lies in
the space of $K_o$-invariant distributions on $M$. It is thus natural
to study $\Delta$ on $K_o$-invariant functions.

Perhaps the most interesting property is the
analytic continuation of the resolvent, which I state before indicating
how it, and other results, relate to many-body scattering.

Fix the branch of the square root function $\sqrt{}$ on 
$\Cx\setminus [0,+\infty)$ which has negative imaginary part when 
$w\in\Cx\setminus[0,+\infty)$. Let $S$ denote that part of the Riemann 
surface for $\lambda \mapsto \sqrt{\lambda-\lambda_0}$ where we continue from 
$\lambda - \lambda_0 \nin [0,+\infty)$ and allow $\arg(\lambda-\lambda_0)$ to change by 
any amount less than $\pi$. In other words, starting in the region 
$\im\sqrt{\lambda-\lambda_0}<0$, we continue across either of the rays where
$\im\sqrt{\lambda-\lambda_0}=0$ and $\re\sqrt{\lambda-\lambda_0}>0$, respectively $< 0$, 
allowing the argument of $\sqrt{\lambda-\lambda_0}$ to change by any amount less
than $\pi/2$ (so that only the positive imaginary axis is not reached). 

\begin{thm}\label{thm:analytic}
With all notation as above, the Green function $G_o(\lambda)$ continues 
meromorphically to $S$ as a distribution. Similarly, as an operator 
between appropriate spaces of $K_o$-invariant functions, the resolvent 
$R(\lambda)$ itself has a meromorphic continuation in this region, with all
poles of finite rank. 
\end{thm}

Having stated the theorem, I indicate how it relates to many-body
scattering.
To do so, fix the point $o$ -- we may as well take it to be the
identity matrix $I$. Now, $M$ is a perfectly nice real analytic
manifold and $\Delta$ is an elliptic operator on it in the
usual sense, so the only question is its behavior at infinity.
In order to describe this, we remark
that any matrix $A\in M$ can be diagonalized,
i.e.\ written as $A=O\Lambda O^t$, with $O\in\SO(3)$ and $\Lambda$ 
diagonal, with entries given by the eigenvalues of $A$.
If $\fraka$ is the set of diagonal matrices of trace $0$,
then $\Lambda \in \exp(\fraka)$. If all eigenvalues of $A$ are distinct,
then $\Lambda$ is determined except for the ordering of the eigenvalues,
and there are only finitely many possibilities for $O$ as well.
However, at the walls, which are defined to be the places where
two eigenvalues coincide, there is much more indeterminacy.
For example, if two eigenvalues coincide, only their joint
eigenspace is well-defined. Correspondingly, we may replace $O$ by
$O'O$ for any $O'\in\SO(3)$ preserving
the eigenspace decomposition and still obtain the desired diagonalization.

This is closely reflected in the structure of the Laplacian at infinity.
In fact, it turns out that on $\SO(3)$-invariant functions,
$\Delta$ is essentially a three-body Hamiltonian on $\fraka$ with first order
interactions and with collision `planes' given by the walls (they
are lines), see e.g.\ \cite[Chapter II, Proposition 3.9]{Helgason:Groups}.
So rather than particles, eigenvalues
scatter in this case! Consequently, many-body results can be adapted
to this setting.

We indicate how this is done. The most succint way of describing the geometry
of $M$ at infinity is to compactify it to a manifold $\bar M$
with codimension two corners.
It has two boundary hypersurfaces, 
$H_\sharp$ and $H^\sharp$, which are perhaps easiest to 
describe in terms of a natural system of local coordinates
derived from the matrix representation of elements in $M$. 
As above, we write $A\in M$ as $A=O\Lambda O^t$,
with $O\in\SO(3)$ and $\Lambda$ 
diagonal.
The ordering of the diagonal entries of $\Lambda$ is 
undetermined, but in the region where no two of them are equal,
we denote them as $0 < \lambda_1 < \lambda_2 < \lambda_3$ (but recall
also that $\lambda_1 \lambda_2 \lambda_3 = 1$). In this region the ratios
\[
\mu = \frac{\lambda_1}{\lambda_2}, \qquad \nu = \frac{\lambda_2}{\lambda_3}
\]
are independent functions, and near the submanifold $\exp(\fraka)$ in $M$ 
we can complete them to a full coordinate system by adding the above-diagonal 
entries $c_{12}$, $c_{13}$, $c_{23}$ in the skew-symmetric matrix $T = \log O$.
On $M$ we have $\mu,\nu>0$, and locally the compactification consists of
replacing $(\mu,\nu)\in(0,1)\times(0,1)$ by $(\mu,\nu)\in[0,1)\times[0,1)$.
Then $H^\sharp = \{\mu=0\}$ and $H_\sharp = \{\nu=0\}$,
and this coordinate system gives the $\Cinf$ structure near the
corner $H_\sharp \cap H^\sharp$.

\begin{figure}[ht]
\begin{center}
\mbox{\epsfig{file=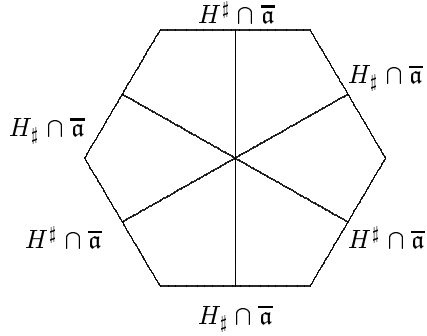}}
\end{center}
\caption{The closure of $\fraka$, or rather $\exp(\fraka)$,
in the compactification
$\bar M$ of $M$. The lines in the interior are the Weyl chamber
walls, playing the
role of collision planes in many-body scattering. The side faces
$H^\sharp\cap\overline{\fraka}$ and $H_\sharp\cap\overline{\fraka}$
correspond to the front faces on Figure~\ref{fig:ressp}. The main
face on Figure~\ref{fig:ressp} would only show up if we did a logarithmic
blow-up of all boundary hypersurfaces of $\bar M$ and then blew up the corner.}
\label{fig:flat1}
\end{figure}

On the other hand, in a neighborhood of 
the interior of $H_\sharp$, for example, we obtain the compactification
and its $\Cinf$ structure as follows.
Write the eigenvalues of $A \in M$, i.e.\ the diagonal entries
of $\Lambda$ in the decomposition for $A$ above, as
$\lambda_1$, $\lambda_2$ and $\lambda_3$. Suppose that $A$ 
lies in a small neighbourhood $\calU$ where 
\[
c< \lambda_1/\lambda_2< c^{-1},\quad \lambda_3>1/c,
\]
for some fixed $c \in (0,1)$. Recall also that $\lambda_3 =
1/\lambda_1 \lambda_2$. These inequalities imply that $\lambda_1=
(\lambda_1/\lambda_2)^{1/2}\lambda_3^{-1/2}<1$ and $\lambda_2= 
(\lambda_2/\lambda_1)^{1/2}\lambda_3^{-1/2}<1$, and $\lambda_3 > 1$
in $\calU$. Hence there is a well-defined decomposition 
$\Real^3 = E_{12} \oplus E_3$ for any $A \in \calU$, 
where $E_{12}$ is the sum of the first two eigenspaces and $E_3$ 
is the eigenspace corresponding to $\lambda_3$, regardless of whether 
or not $\lambda_1$ and $\lambda_2$ coincide. We could write equivalently
$A=O C O^t$, where $C$ is block-diagonal, preserving the splitting
$\Real^2\oplus\Real$
of $\RR^3$. The ambiguity in this factorization is 
that $C$ can be conjugated by an element of $\On(2)$ (acting in the 
upper left block), and $\On(2)$ can be included in the top left
corner of $\SO(3)$ (the bottom right entry being set equal to $\pm 1$
appropriately). Let $C'$ denote the upper-left block of $C$; the
bottom right entry of $C$ is just $\lambda_3$, and so $\lambda_3 \det C' 
= 1$. In other words, $C' = \lambda_3^{-1/2}C''$ where $C''$ is
positive definite and symmetric with determinant $1$, hence
represents an element of $\SL(2)/\SO(2) \equiv \HH^2$. Hence for
an appropriate neighbourhood $\calV$ of $I$ in $\SL(2)/\SO(2)$, 
the neighbourhood $\calU$ is identified with $(\calV \times \SO(3))/\On(2) 
\times (0,c^{3/2})$, where the variable on the last factor is
$s=\lambda_3^{-3/2}$. The compactification then simply
replaces $(0,c^{3/2})_s$ by $[0,c^{3/2})$.
Note that although the action of $\On(2)$ on $\calV$
has a fixed point (namely $I$), its action on $\SO(3)$,
and hence on the product, is free. The neighbourhood $\calV$
can be chosen larger when $\lambda_3$ is larger, and the limiting
`cross-section' $\lambda_3 = \text{const}$ has the form 
$(\HH^2 \times \SO(3))/\mbox{O}(2)$. This space is a fibre bundle
over $\SO(3)/\On(2)$ ($= \RR P^2$) with fibre $\HH^2$. 
Notice that the Weyl chamber wall corresponds to the origin (i.e.\ 
the point fixed by the $\SO(2)$ action) in $\HH^2$. 
I refer to \cite{Mazzeo-Vasy:Sl3} for a more thorough description
of $\bar M$.

On each boundary hypersurface of $M$, it is now easy to describe model
operators for $\Delta$ acting on $\SO(3)$-invariant functions. For instance,
at $H_\sharp$ this model can be considered as an operator $L_\sharp$
on $\Real_s\times\HH^2$, acting on $\SO(2)$-invariant functions.
Explicitly,
\begin{equation*}
L_{\sharp}=\frac{1}{4}(sD_s)^2+i\frac{1}{2}(sD_s)
+\frac{1}{3}\Delta_{\HH^2}.
\end{equation*}
This is tensor product type, so its resolvent can be obtained
from an integral of the resolvents of $\frac{1}{4}(sD_s)^2+i\frac{1}{2}(sD_s)$
and $\frac{1}{3}\Delta_{\HH^2}$. (Note that I am ignoring the weights of
the $L^2$ spaces on which we are working, hence the appearance of the
perhaps strange first
order terms.)

This framework allows one to develop the elliptic theory, for example to
analyze $(\Delta-\lambda)^{-1}$ for
$\lambda\in\Cx\setminus[\lambda_0,+\infty)$. In particular, one
can construct a parametrix for $\Delta$ on $\bar M$ that has a smoothing
error. Since this error has no decay at infinity, it is not compact.
However, the error can be improved by pasting
together the resolvents of $L_\sharp$ and $L^\sharp$,
and applying the result to the error term to remove it modulo a decaying,
hence compact, new error term.
One of the consequences
is then the description of the asymptotic behavior of the Green's function,
see \cite{Mazzeo-Vasy:Sl3}.

The proof of Theorem~\ref{thm:analytic}
relies on the method of complex scaling.
This is an extension of dilation 
analyticity, and was originally developed in the setting of $2$-body 
scattering by Aguilar-Combes \cite{Aguilar-Combes:Analytic} and generalized to
the many-body setting by Balslev-Combes \cite{Balslev-Combes:Spectral}
and later further generalized by G\'erard \cite{Gerard:Distortion}.
We refer to \cite{Hislop-Sigal:Spectral} and \cite[Volume 4]{Reed-Simon}
for an exposition, and to the paper \cite{Sjostrand-Zworski:Complex}
of Sj\"ostrand and Zworski for a slightly different point of view.

The point of complex scaling is to rotate the essential spectrum
of the operator being studied, in this case the Laplacian.
To give the reader a rough idea how this works, consider the hyperbolic
space $\HH^2=\SL(2,\Real)/\SO(2,\Real)$, which may be identified with
the set of two-by-two positive definite matrices $A$ of determinant $1$.
In terms of geodesic normal coordinates $(r,\omega)$ about $o=I$,
the Laplacian is given by
\begin{equation*}
\Delta_{\HH^2}=D_r^2-i(\coth r)D_r+(\sinh r)^{-2}D_\omega^2.
\end{equation*}
Now consider the diffeomorphism $\Phi_\theta:A\mapsto A^w$, $w=e^\theta$,
on $\HH^2$, $\theta\in\Real$. 
This corresponds to dilation along the geodesics through $o$,
since these have the form $\gamma_A:s\mapsto A^{cs}$, $c>0$. Thus,
in geodesic normal coordinates, $\Phi_\theta:
(r,\omega)\mapsto(e^\theta r,\omega)$. $\Phi_\theta$ defines a
group of unitary operators on $L^2(\HH^2)$ via
\begin{equation*}
(U_\theta f)(A)=(\det D_A\Phi_\theta)^{1/2}(\Phi_\theta^*f)(A),
\ J=\det D_A\Phi_\theta=w\frac{\sinh wr}{\sinh r},\ w=e^\theta.
\end{equation*}
Now, for $\theta$ real, consider the scaled Laplacian 
\begin{eqnarray*}
(\Delta_{\HH^2})_\theta & = & U_\theta \Delta_{\HH^2} U_\theta^{-1} =
J^{1/2}\Phi_\theta^* \Delta_{\HH^2} \Phi_{-\theta}^* J^{-1/2} \\
& = & J^{1/2}(w^{-2}D_r^2-iw^{-1}\coth(wr)D_r
+(\sinh(wr))^{-2}D_\omega^2)J^{-1/2}
\end{eqnarray*}
This is an operator on $\HH^2$, with coefficients which extend analytically 
in the strip $|\im\theta|<\frac{\pi}{2}$. The square root is continued 
from the standard branch near $w>0$. 
(The singularity of the coefficients at $r=0$ is only an artifact 
of the polar coordinate representation.) Note
that $(\Delta_{\HH^2})_\theta$ and $(\Delta_{\HH^2})_{\theta'}$ are
unitary equivalent if $\im\theta=\im\theta'$ because of the group properties
of $U_\theta$. The scaled operator, $(\Delta_{\HH^2})_\theta$,
is not elliptic on
all of $\HH^2$ when $0<|\im\theta|<\frac{\pi}{2}$ because for $r$ large
enough, $w^2\sinh(wr)^{-2}$ can lie in $\Real^-$.
However, it is elliptic in some 
uniform neighbourhood of $o$ in $\HH^2$, and its radial part 
\begin{equation*}
(\Delta_{\HH^2})_{\theta,\mathrm{rad}}=J^{1/2}(w^{-2}D_r^2-iw^{-1}\coth(wr)D_r)
J^{-1/2},
\end{equation*}
which corresponds to its action on $\SO(2)$-invariant functions, is
elliptic on 
the entire half-line $r>0$. The model operator for
$(\Delta_{\HH^2})_{\theta,\mathrm{rad}}-\lambda$ at infinity,
\begin{equation*}\begin{split}
&e^{(w-1)r/2}(w^{-2}D_r^2-iw^{-1}D_r-\lambda)e^{-(w-1)r/2}\\
&=e^{(w-1)r/2}((w^{-1}D_r-\frac{i}{2})^2-(\lambda-\frac{1}{4}))e^{-(w-1)r/2},
\end{split}\end{equation*}
is also
invertible on the model space at infinity, $L^2(\Real;e^r\,dr)$, since
this is equivalent to the invertibility of
\begin{equation*}
e^{wr/2}((w^{-1}D_r-\frac{i}{2})^2-(\lambda-\frac{1}{4}))e^{-wr/2}
=w^{-2}D_r^2-(\lambda-\frac{1}{4})
\end{equation*}
on $L^2(\Real;dr)$.
Thus, a parametrix with compact remainder can be 
constructed for $(\Delta_{\HH^2})_{\theta,\mathrm{rad}}$, and this show that 
its essential spectrum lies in $\frac{1}{4}+e^{-2i\im\theta}[0,+\infty)$.
Hence $((\Delta_{\HH^2})_{\theta,\mathrm{rad}}-\lambda)^{-1}$ is
meromorphic outside this set.
In fact, it is well known that there are no poles in this entire strip
(although there are an infinite number on $\arg\sqrt{\lambda-\lambda_0}
= \pi/2$).

Combining this with some more standard technical facts, we are in
a position to apply the theory of Aguilar-Balslev-Combes to prove
that $((\Delta_{\HH^2})_{\mathrm{rad}}-\lambda)^{-1}$, and hence
$(\Delta_{\HH^2}-\lambda)^{-1}$, has an analytic continuation in $\lambda$ 
across $(\frac{1}{4},+\infty)$. This is done by noting that 
for $\SO(2)$-invariant functions $f,g\in L^2(\HH^2)$ and $\theta\in\Real$,
\begin{equation*}
\langle f,((\Delta_{\HH^2})_{\mathrm{rad}}-\lambda)^{-1}g\rangle
=\langle U_{\bar\theta}f,((\Delta_{\HH^2})_{\theta,\mathrm{rad}}-\lambda)^{-1}
U_\theta g\rangle
\end{equation*}
by the unitarity of $U_\theta$. Now if $f$, $g$ lie in a smaller (dense)
class of functions such that $U_\theta f$ and $U_\theta g$ continue
analytically from $\theta\in\Real$, then the meromorphic continuation in $\lambda$ of
the right hand side is obtained by first making $\theta$ complex 
with imaginary part of the appropriate sign, and then allowing $\lambda$ to 
cross the continuous spectrum of $\Delta_{\HH^2}$ without encountering 
the essential spectrum of $(\Delta_{\HH^2})_{\theta,\mathrm{rad}}$.
Hence the left hand side continues meromorphically also. With some additional
care, one can even allow $g$ to be the delta distribution at $o$, yielding the 
meromorphic continuation of the Green's function.

The argument on the higher rank symmetric space $M=\SL(3)/\SO(3)$
is similar. We still use the same scaling
$\Phi_\theta:A\mapsto A^w$, $w=e^\theta$,
on $M$, $\theta\in\Real$. Again, the first concern is that, allowing $\theta$
to become complex, the scaled operator
$\Delta_\theta$ is not elliptic. However, it is elliptic near $o=\Id$,
and the scaled models for it near the walls, such as $(L_\sharp)_\theta$,
remain elliptic at the walls. After all, for the latter, this is just the
ellipticity of $(\Delta_{\HH^2})_\theta$ near the origin, which we have
already observed. This again allows the elliptic parametrix construction
to proceed, supplying the results we needed in order to reach the
framework of complex scaling. This in turn finishes the proof of
Theorem~\ref{thm:analytic}.

\bibliographystyle{plain}
\bibliography{sm}

\end{document}